\documentclass{article}
\def\N{\mathbb N}
\def\Z{\mathbb Z}
\def\R{\mathbb R}
\def\Q{\mathbb Q}

\def\Zmb{{\mathbb Z}_{-\beta}}
\def\A{\mathcal A}
\def\B{\mathcal B}
\def\L{\mathcal L}

\def\pfz{\begin{proof}}
\def\pfk{\end{proof}}

\usepackage{graphicx, epic, amssymb, amsmath,amsthm, color}
\usepackage{multirow}
\newtheorem{thm}{Theorem}[section]
\newtheorem{lem}[thm]{Lemma}
\newtheorem{prop}[thm]{Proposition}
\newtheorem{coro}[thm]{Corollary}
\newtheorem{pozn}[thm]{Remark}

\newtheorem{ex}[thm]{Example}

\begin{document}
\title{Integers in number systems with positive and negative quadratic Pisot base}

\author{Zuzana Mas\'akov\'a and Tom\'a\v s V\'avra\\
Department of Mathematics FNSPE, Czech Technical University in Prague\\
Trojanova 13, 120 00 Praha 2, Czech Republic\\
\tt{emails: zuzana.masakova@fjfi.cvut.cz, t.vavra@seznam.cz}
}
%
\date{\today}
\maketitle
\begin{abstract}
We consider numeration systems with base $\beta$ and $-\beta$, for quadratic Pisot numbers $\beta$ and focus
on comparing the combinatorial structure of the sets $\Z_\beta$ and $\Z_{-\beta}$
of numbers with integer expansion in base $\beta$, resp. $-\beta$.
Our main result is the comparison of languages of infinite words $u_\beta$ and $u_{-\beta}$ coding the ordering of distances between consecutive $\beta$- and $(-\beta)$-integers.
It turns out that for a class of roots $\beta$ of $x^2-mx-m$, the languages coincide, while for other quadratic Pisot numbers the language of $u_\beta$ can be identified only
with the language of a morphic image of $u_{-\beta}$. We also study the group structure of $(-\beta)$-integers.
\end{abstract}
%
%
%

\section*{Introduction}
Numeration systems with negative non-integer base obtained a non-negligible attention since the paper~\cite{ItoSadahiro} of Ito and Sadahiro in 2009.
Some of the articles that followed point out that many properties of these systems are analogous to those of systems with non-integer positive base defined by R\'enyi in 1957~\cite{Renyi}, and are easily derivable. On the other hand, there are interesting results showing that the analogy is sometimes very non-trivial, and that negative base
systems have certain features very different from properties of R\'enyi systems. It is therefore appealing to focus on the comparison of these two types of numeration systems.
The comparison of the dynamical aspects of these systems is the topic of~\cite{kalle}.
The present paper is another attempt in this direction, concentrating on comparing the combinatorial structure of the sets $\Z_\beta$ and $\Z_{-\beta}$
of numbers with integer expansion in base $\beta$, resp. $-\beta$. We perform the study for the case where $\beta$ is a quadratic Pisot number.
For this purpose, we have to put $\beta$- and $(-\beta)$-integers into several different contexts. Thus the first part of the paper can be viewed as a review on properties of
$\Z_\beta$ and $\Z_{-\beta}$ found in~\cite{BuFrGaKr,Thurston,fabre,hedlundmorse40,GuMaPe,ADMP,MaPeVa} and others. In particular, we use the cut-and-project scheme, that describes $\Z_\beta$ and $\Z_{-\beta}$ for $\beta$ being an algebraic unit. We study the infinite words $u_\beta$ and $u_{-\beta}$ coding the ordering of distances in $\Z_\beta$ and $\Z_{-\beta}$ and provide prescriptions for the morphisms/antimorphisms under which these infinite words are invariant.

The main result of the paper is the comparison of languages of $u_\beta$ and $u_{-\beta}$. We show that languages of these infinite words coincide if $\beta$ is a root of $x^2-mx-m$, $m\geq 1$. For all other quadratic Pisot numbers $\beta$ one has to first apply a morphism $\pi$ to the infinite word coding $(-\beta)$-integers in order to obtain equality of the corresponding languages. We use the notion of conjugate morphisms.

In the last part of the paper we concentrate on the arithmetical properties of $\Z_\beta$ and $\Z_{-\beta}$. In particular, we generalize the result of~\cite{ElFrGaVG}, presenting a group operation $\oplus$ on $\Z_\beta$, resp. $\Z_{-\beta}$, with which these sets are isomorphic to $\Z$. The authors of~\cite{ElFrGaVG} define such an operation on $\Z_\beta$ for quadratic Pisot units $\beta$ and show that such an operation is compatible with ordinary addition, i.e.\ whenever the result of $x+y$ is an element of $\Z_{\beta}$, then $x\oplus y=x+y$. By an easy combinatorial argument we show that one can define an operation $\oplus$ compatible with addition on every discrete set with finitely many distances between consecutive points which are linearly independent over $\Q$. We use again the cut-and-project scheme for determining the possible outcomes of $x+y-(x\oplus y)$ for $x,y\in\Z_{-\beta}$
for quadratic Pisot units $\beta$.

The organization of the paper is as follows. In Section~\ref{sec:preli} we recall necessary notions from combinatorics on words, present the numeration systems with positive and negative base and define the $\beta$- and $(-\beta)$-integers. In Section~\ref{sec:cap} we show the simplest cut-and-project scheme that can be applied for studying $\beta$- and $(-\beta)$-integers in the quadratic Pisot case. We cite the results identifying $\Z_\beta$ and $\Z_{-\beta}$ as cut-and-project sets. Section~\ref{sec:combi} introduces the infinite words $u_\beta$ and $u_{-\beta}$ coding the distances between consecutive $\beta$- and $(-\beta)$-integers. These words are invariant under morphisms, resp.\ antimorphisms. We list the values of the distances and provide the prescriptions of the corresponding  morphisms, resp. antimorphisms, by citing the known results and giving proofs where such are nowhere explicitly given. The main results of this paper are in Sections~\ref{sec:vztah} and~\ref{sec:group}. Section~\ref{sec:vztah} demonstrates the geometric and combinatorial relation of
$\Z_\beta$ and $\Z_{-\beta}$. We compare the languages of the infinite words $u_\beta$ and $u_{-\beta}$. Section~\ref{sec:group} concentrates on arithmetical aspects of $(-\beta)$-integers.


\section{Preliminaries}\label{sec:preli}

An alphabet $\A$ is a finite set of symbols, called letters. By $\A^*$ we denote the monoid of all finite words over the alphabet $\A$ equipped
with the empty word $\epsilon$ and the operation of concatenation. The concatenation of $k$ copies of a word $w$ is denoted by $w^k = w\cdots w$,
for $k\in\N$. The length of a finite word $w=w_1\cdots w_n$, $w_i\in\A$, is denoted by $|w|=n$.
By $|w|_a$ we denote the number of occurrences of a letter $a$ in the word $w$.
We will consider one-directional infinite words $u=u_0u_1u_2\cdots$, resp. $u=\cdots u_{-2}u_{-1}$, but also pointed bidirectional infinite words
$u=\cdots u_{-2}u_{-1}|u_0u_1u_2\cdots$. The delimiter $\mid$ marks the position of the origin.
A word $w$ is a factor of a word $u$, if $u=vwv'$ for some (finite or infinite) words $v,v'$. If $v=\epsilon$, then $w$ is a prefix of $u$, if
$v'=\epsilon$, then $w$ is a suffix of $u$.
The language of an infinite word $u$ is the set ${\mathcal L}(u)$ of all factors of $u$. The factor complexity of the word $u$ is the function ${\mathcal C}:\N\to\N$, where ${\mathcal C}(n)$ is defined as the number of factors of $u$ of length $n$.

The infinite repetition of a finite word $w$ is denoted by $w^\omega=www\cdots$.
A one-directional infinite word $u$ is said eventually periodic if $u=vw^\omega=vwww\cdots$, for some finite words $v,w$, otherwise it is called aperiodic.
It is not difficult to show that an aperiodic one-directional word must satisfy ${\mathcal C}(n)\geq n+1$ for $n\in\N$. Aperiodic one-directional infinite words
of minimal complexity ${\mathcal C}(n)= n+1$ for $n\in\N$ are called Sturmian words. Sturmian words are by definition binary words (for ${\mathcal C}(1)=2$) and have many equivalent definitions, for a recent overview of them see~\cite{jungle}. We will be interested in the following one. A one-directional infinite word over an alphabet $\A$ is called balanced if $\big||w|_a-|v|_a\big|\leq 1$ for every pair of factors $w,v\in{\mathcal L}(u)$ of the same length $|w|=|v|$ and every $a\in\mathcal A$. A binary aperiodic word is Sturmian if and only if it is balanced, see~\cite{hedlundmorse40}.

The definition of Sturmian words can be extended to bidirectional infinite words, by requiring aperiodicity and balance, see~\cite{pytheas}. The language of such a Sturmian bidirectional word $u$ is the same as the language of any of its infinite suffix, which is a Sturmian word in the original sense.

 For finite alphabets $\A$ and $\B$, a morphism $\sigma:\A^*\to\B^*$ satisfies $\sigma(vw)=\sigma(v)\sigma(w)$ for any $v,w\in\A^*$. Its action can be extended to both one-directional and pointed bidirectional infinite words by
$$
\begin{aligned}
\sigma(u_0u_1u_2\cdots) &= \sigma(u_0)\sigma(u_1)\sigma(u_2)\cdots\,, \\
\sigma(\cdots u_{-2}u_{-1}|u_0u_1u_2\cdots) &= \cdots \sigma(u_{-2})\sigma(u_{-1})|\sigma(u_0)\sigma(u_1)\sigma(u_2)\cdots\,.
\end{aligned}
$$
If the morphism $\sigma:\A^*\to\A^*$ is non-erasing, i.e.\ $\sigma(a)\neq \epsilon$ for all $a\in\A$, and there exists a letter $a\in\A$ such that $\sigma(a)=aw$ for some non-empty word $w$, then it is called a substitution. A substitution has always a fixed point, i.e.\ one-directional infinite word $u$ such that $\sigma(u)=u$, namely $u=\lim_{j\to\infty}\sigma^j(a)$. Existence of a bidirectional fixed point is ensured requiring moreover a letter $b$ such that $\sigma(b)=vb$
for some non-empty word $v$. Then the fixed point is $u=\lim_{n\to\infty}\sigma^n(b)|\sigma^n(a)$. The limit is taken over the product topology on $\A^\N$, resp. $\A^\Z$.

In a similar way, one defines antimorphisms on $\A^*$ by $\overleftarrow{\sigma}(vw)=\overleftarrow{\sigma}(w)\overleftarrow{\sigma}(v)$ for any $v,w\in\A^*$. A fixed point of an antimorphism $\overleftarrow{\sigma}$ is a bidirectional infinite word $\cdots u_{-2}u_{-1}|u_0u_1u_2\cdots$ such that
$$
\overleftarrow{\sigma}(\cdots u_{-2}u_{-1}|u_0u_1u_2\cdots) = \cdots \overleftarrow{\sigma}(u_2)\overleftarrow{\sigma}(u_1)\overleftarrow{\sigma}(u_0)|\overleftarrow{\sigma}(-1)\overleftarrow{\sigma}(-2)\cdots\,.
$$
Note that a second iteration of an antimorphism is a morphism.

Given a morphism $\sigma$ over $\A=\{a_1,\dots,a_d\}$, one defines its incidence matrix $M_\sigma$, $(M_\sigma)_{ij} = |\sigma(a_j)|_{a_i}$. A~morphism is called primitive, if there is a constant $k\in\N$ such that $M_\sigma^k>0$. The eigenvalues of the matrix $M_\sigma$ predicate for example about balance properties of the fixed point of $\sigma$.

\bigskip
We will be interested in properties of infinite words associated with the set of numbers with integer expansion in an irrational base.
Let $\alpha\in\R$, $|\alpha|>1$, and let $\A$ be a finite subset of $\R$. An $(\alpha,\A)$-representation
of a real number $x$ is the series
$$
x=\sum_{i\leq k} x_i\alpha^i\,,\qquad x_i\in\A\,.
$$
Sometimes we may write $x$ equal to the sequence of digits $x_i$, namely
$$
x=x_kx_{k-1}\cdots x_1x_0\text{\small$\bullet$}x_{-1}x_{-2}\cdots\,.
$$

A specific $(\alpha,\A)$-representation with the set of digits $\A\subset\Z$ is obtained as follows.
Consider the transformation $T:I\to I$ on an interval $I=[l,l+1)$, given by $T(x)=\alpha x - \lfloor\alpha x - l\rfloor$. For an $x\in I$ set
$$
x_{-i}= \lfloor \alpha T^{i-1}(x) - l\rfloor\,,\qquad \text{ for $i\geq 1$}\,.
$$
Then
$$
x=\frac{x_{-1}}{\alpha}+\frac{x_{-2}}{\alpha^2}+\frac{x_{-3}}{\alpha^3}+\cdots.
$$
The string of digits is then denoted by $d(x)=x_{-1}x_{-2}x_{-3}\cdots$.

A string $y_1y_2y_3\cdots$ of elements of the digit set $\A$ is called $(\alpha,\A)$-admissible, if it is equal to the  $(\alpha,\A)$-representation $d(x)$ for some $x\in I$. One can derive for example from~\cite{Gora} that a string $y_1y_2y_3\cdots$ of integers is $(\alpha,\A)$-admissible if and only if every suffix $y_iy_{i+1}y_{i+2}\cdots$ satisfies
\begin{equation}\label{eq:admisobec}
d(l) \preceq y_iy_{i+1}y_{i+2}\cdots \prec \lim_{\delta \to 0+} d(l+1-\delta)\,,
\end{equation}
where $\preceq$ stands for the lexicographical order $\preceq_{\text{\tiny lex}}$ when $\alpha>0$ and for the alternate order $\preceq_{\text{\tiny alt}}$  when $\alpha<0$.
Recall that $x_1x_2x_3 \cdots \preceq_{\text{\tiny alt}} y_1y_2y_3 \cdots$ if $(-1)^k x_k < (-1)^k y_k$, where $k =  \min\{i\mid x_i \neq y_i\}$.

It is not difficult to show that the corresponding order $\preceq$ reflects the natural order on reals in the interval $I$, namely that
$$
d(x)\preceq d(y) \quad\iff\quad x\leq y\,.
$$

The two specific choices of numeration systems which we consider here are the R\'enyi system (cf.~\cite{Renyi}) and the Ito-Sadahiro system (cf.~\cite{ItoSadahiro}) which share many properties but display also important differences.


In the R\'enyi system we set $\alpha=\beta>1$ and $I=[0,1)$. The digits $x_i$ then take values in the set $\A=\A_\beta=\{d\in\Z \mid 0\leq d <\beta\}$. For the sake of clearness, we denote $d(x)=d_\beta(x)$ for $x\in[0,1)$.
The admissibility condition~\eqref{eq:admisobec} states that a string of integers $y_1y_2y_3\cdots$ is admissible in the R\'enyi system if and only if for every $i\geq 1$, we have
\begin{equation}\label{eq:renyiadmiss}
0^\omega\preceq_{\text{\tiny lex}} y_iy_{i+1}y_{i+2}\cdots \prec_{\text{\tiny lex}} \lim_{\delta \to 0+} d_\beta(1-\delta)\,.
\end{equation}
This condition was first derived by Parry~\cite{Parry}.

In order to expand any positive number $x\in\R$ in the R\'enyi system, one uses the expansion of $x/\beta^k,$ where $k\in\mathbb Z$ such that $x/\beta^k\in[0,1).$ Then we have
$$
d(x/\beta^k)=x_{-1}x_{-2}x_{-3}\dots\quad\Rightarrow\quad x=x_{-1}x_{-2}\dots x_{-k}\bullet x_{-k-1}\dots.
$$
This representation of $x$ is called the $\beta$-expansion and is denoted $\langle x\rangle_{\beta}$.
We focus on the set of numbers with integer $\beta$-expansion, namely the set
$$
\Z_{\beta}=\Z_\beta^+\cup-\Z_\beta^+\,,\qquad\text{where}\quad \Z_\beta^+=\{ x\geq 0 \mid \langle x\rangle_{\beta} = x_k\cdots x_{1}x_0\bullet 0^\omega\}\,.
$$

The Ito-Sadahiro number system allows one to expand every real number $x$ (even negative) without the use of a sign.
We set $\alpha=-\beta<-1$, $I=[l,l+1)$, where $l=-\beta/(\beta+1)$, i.e.\ $I=\big[-\beta/(\beta+1),1/(\beta+1)\big)$. The digit set is now equal to $\A=\A_{-\beta}=\{d\in\Z \mid 0\leq d \leq \beta\}$. Here we denote $d(x)=d_{-\beta}(x)$ for $x\in[l,l+1)$. A string of integers $y_1y_2y_3\cdots$ is admissible in the Ito-Sadahiro system if and only if for every $i\geq 1$, we have
\begin{equation}\label{eq:ISadmiss}
d_{-\beta}(l)\preceq_{\text{\tiny alt}} y_iy_{i+1}y_{i+2}\cdots \prec_{\text{\tiny alt}} \lim_{\delta \to 0+} d_\beta(l+1-\delta)\,.
\end{equation}
In~\cite{ItoSadahiro} it is shown that $\lim_{\delta \to 0+} d_\beta(l+1-\delta)$ is strongly dependent on the string $d(l)$, namely $\lim_{\delta \to 0+} d_\beta(l+1-\delta) = \big(0d_1\cdots d_{q-1}(d_{q}-1)\big)^\omega$ if $d(l)=(d_1\cdots d_{q-1} d_q)^\omega$ is purely periodic with odd period-length $q$, and $\lim_{\delta \to 0+} d_\beta(l+1-\delta) = 0d(l)$ otherwise.

In order to expand any real number $x\in \R$, take a suitable power of $(-\beta)$ such that $x/(-\beta)^k\in I^\circ=\big(-\beta/(\beta+1),1/(\beta+1)\big)$. If $d(x/(-\beta)^k)=x_{-1}x_{-2}x_{-3}\dots$ then we put
$\langle x\rangle_{-\beta}=x_{-1}\dots x_{-k}\bullet x_{-k-1}\dots$.
The set of $(-\beta)$-integers is defined by
$$
\Z_{-\beta}=\{ x\in\R \mid \langle x\rangle_{-\beta} = x_k\cdots x_{1}x_0\bullet 0^\omega\}\,.
$$

One can determine the distances between consecutive elements of $\Z_\beta$ and $\Z_{-\beta}$ and, in case that the distances take only finitely many values, one can code their ordering by an infinite word over a finite alphabet. It is interesting to study combinatorial properties of these infinite words, such as invariance under morphism. We will describe these infinite words for quadratic Pisot numbers $\beta$ and use them to relate $\Z_\beta$ and $\Z_{-\beta}$.

\section{${\beta}$- and ${(-\beta)}$-integers as cut-and-project sequences}\label{sec:cap}

In this paper we focus on the class of quadratic Pisot numbers. A Pisot number $\beta$ is the root $>1$ of a monic irreducible polynomial with integer coefficients whose other roots (the algebraic conjugates of $\beta$) are in modulus smaller than 1.
In this paper we focus on the quadratic case.
It can be easily shown that quadratic Pisot numbers are precisely the larger roots of the polynomials
\begin{equation}\label{eq:kvadrpisot}
\begin{aligned}
x^2-mx-n\,,\qquad&\ m\geq n\geq 1\,,\\
x^2-mx+n\,,\qquad&\ m-2\geq n\geq 1\,.
\end{aligned}
\end{equation}
A special role is played by quadratic Pisot units, i.e.\ roots of~\eqref{eq:kvadrpisot} with $n=1$.

For any algebraic number $\beta$, the minimal subfield of ${\mathbb C}$ containing $\beta$ is denoted by $\Q(\beta)$. For quadratic
$\beta$ we have $\Q(\beta)=\{a+b\beta\mid a,b\in\Q\}$. The only non-trivial automorphism over such a field is the Galois automorphism
$':x=a+b\beta\mapsto x'=a+b\beta'$, where $\beta'$ is the conjugate of $\beta$.

In many cases, it is advantageous to apply the Galois conjugation in the study of numeration systems. In particular,
we will find useful to study $\{x'\mid x\in\Z_\beta\}$, resp. $\{x'\mid x\in\Z_{-\beta}\}$.
Obviously, if $\beta$ is a Pisot number, this set is bounded as $\beta$-integers have in their $\beta$-expansion only non-negative powers of $\beta$. The same holds for $(-\beta)$-integers. However, if $\beta$ is a Pisot unit, one can say more, using the
so-called cut-and-project scheme.

For given irrational numbers $\varepsilon, \eta$, $\varepsilon\neq\eta$, consider vectors $\vec{x}_1=\frac1{\varepsilon-\eta}(\varepsilon,-1)$,  $\vec{x}_2=\frac1{\eta-\varepsilon}(\eta,-1)$. Then, $(a+b\eta)\vec{x}_1+(a+b\varepsilon)\vec{x}_2=(a,b)$.
Consequently, for a given point $(a,b)\in\Z^2$, the value $\pi_1(a,b)=a+b\eta$ is its projection to the line $V_1=\R\vec{x}_1$;
the value $\pi_2(a,b)=a+b\varepsilon$ is its projection to the line $V_2=\R\vec{x}_2$.
Since  $\varepsilon, \eta$ are irrational, the mapping $\pi_1$ restricted to $\Z^2$ is an injection, and we can define $^*:\Z+\Z\eta \to
\Z+\Z\varepsilon$ by $^*=\pi_2\circ\pi_1^{-1}$.
We have thus the cut-and-project scheme,
$$
\begin{array}{ccccc}
V_1 &\stackrel{\pi_1}{\longleftarrow}& \R^2=V_1\times V_2 &\stackrel{\pi_2}{\longrightarrow}& V_2\\
&&\cup&&\\
\Z+\Z\eta&\stackrel{\pi_1}{\longleftarrow}& \Z^2 &\stackrel{\pi_2}{\longrightarrow}& \Z+\Z\varepsilon
\end{array}
$$
For a bounded set $\Omega\subset\R$ we define a cut-and-project set $\Sigma_{\eta,\varepsilon}(\Omega)$ by
$$
\Sigma_{\varepsilon,\eta}(\Omega) = \{x\in\Z+\Z\eta \mid x^*\in\Omega\}\,.
$$
For our purposes, it suffices to consider $\Omega$ to be an interval.

In general, a cut-and-project set can be defined by projection of certain lattice points $z\in L\subset\R^{r+s}=V_1\times V_2$ to a suitably oriented $r$-dimensional subspace $V_1$, where the choice of points to be projected is directed by the projection to the subspace $V_2$. For details and a list of references, see~\cite{Moody}.
For our purposes, it is sufficient to limit our considerations to projection of $\Z^2$ to one-dimensional subspaces. Such cut-and-project schemes are subject of~\cite{GuMaPe}. One of the results shown there is that the distances between consecutive points of $\Sigma_{\varepsilon,\eta}(\Omega)$ take two or three values and that the corresponding infinite word is a coding of exchange of two or three intervals. In particular, we will use the following statement which can be derived from~\cite{GuMaPe}.

\begin{thm}\label{t:zgumape}
Let $\varepsilon, \eta$ be distinct irrational numbers and let $\Omega$ be a non-degenerated interval containing $0$. If
$\Sigma_{\varepsilon,\eta}(\Omega)=\{t_j\mid j\in\Z\}\subset\R$, $t_j<t_{j+1}$, $t_0=0$, then there exist positive values $\Delta_0$, $\Delta_1$ such that the distances between consecutive elements of $\Sigma_{\varepsilon,\eta}(\Omega)$ take values $t_{j+1}-t_j\in\{\Delta_0,\Delta_1, \Delta_0+\Delta_1\}$. Moreover, if $t_{j+1}-t_j\in\{\Delta_0,\Delta_1\}$, then the pointed bidirectional infinite word $u_{\varepsilon,\eta}(\Omega)=\cdots u_{-2}u_{-1}|u_0u_1u_2\cdots$ over $\{0,1\}$ defined by
$u_j=X$ if $t_{j+1}-t_j=\Delta_X$ is a bidirectional Sturmian word.
\end{thm}

In connection to $\beta$- and $(-\beta)$-integers we will need a special case of the set $\Sigma_{\varepsilon,\eta}(\Omega)$, namely
for $\varepsilon=\beta'$ and $\eta=\beta$ for a given quadratic Pisot unit $\beta$. Then $\Z+\Z\eta=\Z[\beta]$ and $\Z+\Z\varepsilon=\Z[\beta']=\Z[\beta]$, and the mapping $^*$ coincides with the Galois conjugation.
For a bounded set $\Omega\subset\R$ we thus have the cut-and-project set
$$
\Sigma_\beta(\Omega) = \Sigma_{\beta',\beta}(\Omega) = \{x\in\Z[\beta] \mid x'\in\Omega\}\,,
$$
where we simplify the notation in indices. We will need the following properties of $\Sigma_\beta(\Omega)$. The first of the statements of Proposition~\ref{l:cap} would be trivially valid for any cut-and-project set, the last two are specific for the algebraic choice of
parameters $\varepsilon$ and $\eta$.

\begin{prop}\label{l:cap}
Let $\alpha,x_0\in\Z[\beta]$ and let $\alpha$ be an algebraic unit. Let $\Omega,\Omega_1,\Omega_2$ be intervals such that $\Omega=\Omega_1\cup\Omega_2$. Then
\begin{enumerate}
\item $\Sigma_\beta(\Omega_1)\cup\Sigma_\beta(\Omega_2)=\Sigma_\beta(\Omega_1\cup\Omega_2)$,
\item $x_0+\Sigma_\beta(\Omega) = \Sigma_\beta(x'_0+\Omega)$,
\item $\Sigma_\beta(\alpha'\Omega)=\alpha\Sigma_\beta(\Omega)$.
\end{enumerate}
\end{prop}

\pfz
The first statement follows directly from the definition. For the second statement, we have
$$
x_0+\Sigma_\beta(\Omega) = \{x_0+x\in\Z[\beta]\ |\ x'\in\Omega\} =\{x\in\Z[\beta]\ |\ (x-x_0)'\in\Omega\} = \Sigma_\beta(x'_0+\Omega)\,.
$$
For the third one, we have
$$
\Sigma_\beta(\alpha'\Omega)=\{x\in\Z[\beta]\ |\ x'\in\alpha'\Omega\}=
\{x\in\Z[\beta]\ |\ (\alpha^{-1}x)'\in\Omega\}=\{\alpha x\in\Z[\beta]\ |\ x'\in\Omega\}=\alpha\Sigma_\beta(\Omega),
$$
where we have used that $\alpha$ is an algebraic unit and hence $\alpha\Z[\beta]=\Z[\beta]$. This means that $\alpha x\in \Z[\beta]$
is equivalent to $x\in\Z[\beta]$.
\pfk

\begin{figure}[ht]
{\setlength{\unitlength}{0.9pt}
\qquad\qquad
\begin{picture}(200,255)
\put(0,10){\line(1,0){180}}
\put(0,30){\line(1,0){180}}
\put(180,102){$a$}
\put(0,50){\line(1,0){180}}
\put(0,70){\line(1,0){180}}
\put(0,90){\line(1,0){180}}
\put(70,110){\vector(1,0){120}}
\put(70,110){\vector(-1,0){80}}
\put(0,130){\line(1,0){180}}
\put(0,150){\line(1,0){180}}
\put(0,170){\line(1,0){180}}
\put(0,190){\line(1,0){180}}
\put(0,210){\line(1,0){180}}
\put(0,230){\line(1,0){180}}
\put(10,0){\line(0,1){240}}
\put(60,240){$b$}
\put(30,0){\line(0,1){240}}
\put(50,0){\line(0,1){240}}
\put(70,110){\vector(0,1){140}}
\put(70,110){\vector(0,-1){120}}
\put(90,0){\line(0,1){240}}
\put(110,0){\line(0,1){240}}
\put(130,0){\line(0,1){240}}
\put(150,0){\line(0,1){240}}
\put(170,0){\line(0,1){240}}
{\put(10,10){\circle*{5}}}
{\put(10,30){\circle*{5}}}
{\put(30,30){\circle*{5}}}
{\put(10,50){\circle*{5}}}
{\put(30,50){\circle*{5}}}
{\put(50,50){\circle*{5}}}
{\put(30,70){\circle*{5}}}
{\put(50,70){\circle*{5}}}
{\put(30,90){\circle*{5}}
{\put(50,90){\circle*{5}}}
{\put(70,90){\circle*{5}}}
{\put(50,110){\circle*{5}}}
{\put(70,110){\circle{5}}}
{\put(90,110){\circle*{5}}}
{\put(70,130){\circle*{5}}}
{\put(90,130){\circle*{5}}}
{\put(90,150){\circle*{5}}}
{\put(110,150){\circle*{5}}}}
{\put(110,130){\circle*{5}}}
{\put(90,170){\circle*{5}}}
{\put(110,170){\circle*{5}}}
{\put(130,170){\circle*{5}}}
{\put(110,190){\circle*{5}}}
{\put(130,190){\circle*{5}}}
{\put(150,190){\circle*{5}}}
{\put(130,210){\circle*{5}}}
{\put(150,210){\circle*{5}}}
{\put(130,230){\circle*{5}}}
{\put(150,230){\circle*{5}}}
{\put(170,230){\circle*{5}}}
%
\drawline(10,146)(150,62)
\drawline(54.5,119)(131,241)
\drawline(94,96)(182,236)
\drawline(46,124)(0,49)
\drawline(84,102)(22.5,0)
\end{picture}
\hfill
\begin{picture}(200,255)
\put(0,10){\line(1,0){180}}
\put(0,30){\line(1,0){180}}
\put(180,102){$a$}
\put(0,50){\line(1,0){180}}
\put(0,70){\line(1,0){180}}
\put(0,90){\line(1,0){180}}
\put(70,110){\vector(1,0){120}}
\put(70,110){\vector(-1,0){80}}
\put(0,130){\line(1,0){180}}
\put(0,150){\line(1,0){180}}
\put(0,170){\line(1,0){180}}
\put(0,190){\line(1,0){180}}
\put(0,210){\line(1,0){180}}
\put(0,230){\line(1,0){180}}
\put(10,0){\line(0,1){240}}
\put(40,240){$b$}
\put(30,0){\line(0,1){240}}
\put(70,0){\line(0,1){240}}
\put(50,110){\vector(0,1){140}}
\put(50,110){\vector(0,-1){120}}
\put(90,0){\line(0,1){240}}
\put(110,0){\line(0,1){240}}
\put(130,0){\line(0,1){240}}
\put(150,0){\line(0,1){240}}
\put(170,0){\line(0,1){240}}
{\put(30,10){\circle*{5}}}
{\put(10,10){\circle*{5}}}
{\put(10,30){\circle*{5}}}
{\put(30,30){\circle*{5}}}
{\put(30,50){\circle*{5}}}
{\put(50,50){\circle*{5}}}
{\put(30,70){\circle*{5}}}
{\put(50,70){\circle*{5}}}
{\put(70,70){\circle*{5}}
{\put(50,90){\circle*{5}}}
{\put(70,90){\circle*{5}}}
{\put(50,110){\circle{5}}}
{\put(70,110){\circle*{5}}}
{\put(90,110){\circle*{5}}}
{\put(70,130){\circle*{5}}}
{\put(90,130){\circle*{5}}}
{\put(90,150){\circle*{5}}}
{\put(110,150){\circle*{5}}}}
{\put(110,130){\circle*{5}}}
{\put(90,170){\circle*{5}}}
{\put(110,170){\circle*{5}}}
{\put(130,170){\circle*{5}}}
{\put(110,190){\circle*{5}}}
{\put(130,190){\circle*{5}}}
{\put(150,190){\circle*{5}}}
{\put(130,210){\circle*{5}}}
{\put(150,210){\circle*{5}}}
{\put(130,230){\circle*{5}}}
{\put(150,230){\circle*{5}}}
{\put(170,230){\circle*{5}}}
\drawline(10,146)(150,62)
\drawline(0,29)(131,241)
\drawline(34.5,0)(182,236)
%
\end{picture}
\qquad\qquad
\caption{The cut-and-project scheme for the set of $\tau$-integers (left) and $(-\tau)$-integers (right). Every element $x\in \Z_\beta$ and $y\in\Z_{-\beta}$ can be written in the form $x=a+b\tau$, $y=c+d\tau$, where $a,b,c,d\in\Z$. The picture represents $x$ as a pair $(a,b)\in\Z^2$, resp. $(c,d)\in\Z^2$.}
\label{f:c&p}
}
\end{figure}
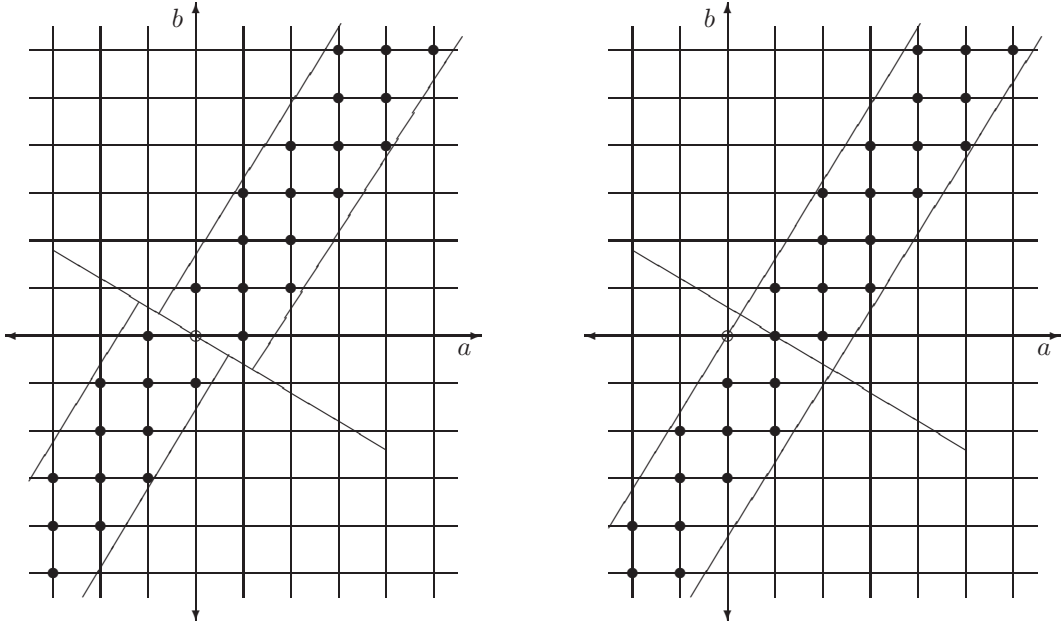

As shown in~\cite{BuFrGaKr}, for quadratic Pisot units $\beta$, the positive part $\Z_\beta^+$ of $\beta$-integers can be identified with the positive part of a cut-and-project set. The reason why one cannot write the equality in the following proposition for the entire set $\Z_\beta$ is that $\Z_\beta$ is in some sense `artificially' defined so that it is symmetric with respect to zero.

\begin{prop}\label{p:positCP1}
Let $\beta>1$ be the root of $x^2-mx-1$, $m\geq 1$. Then
$$
\Z_\beta^+ = \Sigma_\beta\big((-1,\beta)\big)\cap \mathbb R^+ = \big\{x\in\Z[\beta] \,\big|\, x'\in (-1,\beta) \big\} \cap \mathbb R^+\,.
$$
Let $\beta>1$ be the root of $x^2-mx+1$, $m\geq3$. Then
$$
\Z_\beta^+ = \Sigma_\beta\big([0,\beta)\big)\cap \mathbb R^+ =\big\{x\in\Z[\beta] \,\big|\, x'\in [0,\beta) \big\} \cap \mathbb R^+\,.
$$
\end{prop}

Unlike the case of positive base, in the Ito-Sadahiro system, one can represent negative numbers without the use of a sign.
This allows to identify the entire set $\Z_{-\beta}$ as a cut-and-project set. The difference between integers in positive and negative base is illustrated in Figure~\ref{f:c&p} for the case $\beta=\tau=\frac12(1+\sqrt5)$, the golden ratio.
The following proposition appeared as Lemma~5 and~9 in~\cite{MaPePurPer}.

\begin{prop}\label{p:negatCP1}
Let $\beta>1$ be the root of $x^2-mx-1$, $m\geq 1$. Then
$$
\Z_{-\beta} =
\begin{cases}
\Sigma_\beta\big([0,\beta)\big)=\big\{z\in\Z[\beta] \mid z' \in [0, \beta)\big\} & \text{ for }m\geq 2 \,,\\
\Sigma_\beta\big([0,\beta^2)\big)=\big\{z\in\Z[\beta] \mid z' \in [0, \beta^2)\big\} & \text{ for }m= 1 \,.
\end{cases}
$$
Let $\beta>1$ be the root of $x^2-mx+1$, $m\geq3$. Then
$$
\Z_{-\beta} = \Sigma_\beta\Big(\big(-\tfrac{\beta-1}{\beta+1},\beta\tfrac{\beta-1}{\beta+1}\big)\Big)= \Big\{z\in\Z[\beta] \,\Big|\, z' \in \big(-\tfrac{\beta-1}{\beta+1},\beta\tfrac{\beta-1}{\beta+1}\big)\Big\}\,.
$$
\end{prop}

\section{Combinatorial properties of $(-\beta)$-integers}\label{sec:combi}

Before starting to compare $\Z_\beta$, resp. $\Z_{-\beta}$ we shall state the necessary properties of $\beta$- and $(-\beta)$-integers. We cite the known ones and provide proofs for the rest.

Recall that the distances between consecutive points of $\Z_\beta$ can be described using the infinite word $\lim_{\delta\to0+}d_\beta(1-\delta)$ of~\eqref{eq:renyiadmiss}, see~\cite{Thurston}. They take finitely many values if and only if this sequence is eventually periodic. Bases $\beta$, for which this happens are called Parry numbers. Similarly, the distances between consecutive points of $\Z_{-\beta}$ can be described using the infinite word $d_{-\beta}(l)$ of~\eqref{eq:ISadmiss}, cf.~\cite{ADMP}. They take finitely many values for the Ito-Sadahiro numbers $\beta$ (in~\cite{steinerSubstRairo} called Yrrap numbers), i.e.\ such that $d_{-\beta}(l)$ is eventually periodic.

Although the class of Parry numbers and the class of Ito-Sadahiro numbers do not coincide, as shown in~\cite{Lingping}, Pisot numbers belong to both classes. In fact, assuming that $\beta$ is a quadratic number, the notions of Parry numbers, Ito-Sadahiro numbers and Pisot numbers coincide~\cite{Bassino,MaPeVa}.
If $\beta$ is a Pisot number, the distances between consecutive points of $\Z_\beta$, resp. $\Z_{-\beta}$ take only finitely many values. In particular,
for quadratic Pisot numbers, the distances are only two.

\begin{prop}\label{p:mezeryposit}
Let $\beta>1$ be a quadratic Pisot number. Let $\Z_\beta=\{t_j\mid j\in\Z\}$, where $t_0=0$ and $t_j<t_{j+1}$.
Then the distances between consecutive $\beta$-integers take values $t_{j+1}-t_j \in \{\Delta_0^+,\Delta_1^+\}$, where
$$
\Delta_0^+=1 \quad \text{ and }\quad
\Delta_1^+=\begin{cases}
\frac{n}\beta & \text{ if } \beta^2=m\beta+n,\ m\geq n\geq 1\,,\\
1-\frac{n}\beta & \text{ if } \beta^2=m\beta-n,\ m-2\geq n\geq 1\,.\\
\end{cases}
$$
\end{prop}

The ordering of distances in $\Z_\beta$ can be coded by an infinite word $u_\beta$. Since $\Z_\beta$ is symmetric with respect to 0, we define $u_\beta$ as a one-sided
infinite word
\begin{equation}\label{eq:defuposit}
u_{\beta}=u_{0}u_1u_2\cdots\in\{0,1\}^\N\,,\qquad \text{where \ $u_j=X$ \ if \ $t_{j+1}-t_j=\Delta_X^+$.}
\end{equation}
It is known that the infinite word $u_\beta$ is a fixed point of the canonical substitution~\cite{fabre}.

\begin{prop}\label{p:substposit}
Let $\beta>1$ be a quadratic Pisot number. Let $\varphi_\beta:\{0,1\}^*\to\{0,1\}^*$ be the canonical morphism
$$
\begin{array}{llll}
\varphi_\beta(0)=0^{m}1\,, &\varphi_\beta(1)=0^{n}\,,\qquad &\text{ if } \beta^2=m\beta+n,\ &m\geq n\geq 1\,,\\
\varphi_\beta(0)=0^{m-1}1\,, &\varphi_\beta(1)=0^{m-n-1}1\,,\qquad &\text{ if } \beta^2=m\beta-n,\ &m-2\geq n\geq 1\,.
\end{array}
$$
Then the one-sided infinite word $u_\beta$ coding the sequence of non-negative $\beta$-integers $\Z_\beta^+$ is a fixed point of $\varphi_\beta$.
In particular, we have $u_{\beta}=\lim_{j\to +\infty}\varphi_\beta^j(0)$.
\end{prop}

An infinite word $u_\beta$ coding positive $\beta$-integers and the corresponding canonical morphism $\varphi_\beta$ fixing $u_\beta$ can be defined for every Parry number
$\beta$. Among all such words, exactly those corresponding to quadratic Pisot units are Sturmian.

\begin{prop}\label{p:sturmianposit}
Let $u_{\beta}$ be the one-directional infinite word coding $\Z_{\beta}^+$. Then $u_\beta$ is Sturmian
if and only if $\beta$ is a quadratic Pisot unit.
\end{prop}

Statements similar to Propositions~\ref{p:mezeryposit}, \ref{p:substposit} and~\ref{p:sturmianposit} can be stated for $(-\beta)$-integers.
Before presenting these results, recall that the admissibility of digit strings as $(-\beta)$-expansions is decided by the alternate order condition~\eqref{eq:ISadmiss} using the expansion $d_{-\beta}(l)$ of the left-end point $l=-\frac\beta{\beta+1}$ of the interval $I=[l,l+1)$. In case of quadratic Pisot numbers these strings are given as follows,
\begin{align}
d_{-\beta}(l) &= m(m-n)^\omega\,,  &\text{ for }\beta^2=m\beta+n,\ m\geq n\geq 1\,,\quad\ \ \label{eq:rozvojela}\\
d_{-\beta}(l) &= \big((m-1)n\big)^\omega\,,  &\text{ for }\beta^2=m\beta-n,\ m-2\geq n\geq 1\,.\label{eq:rozvojelb}
\end{align}

For the class~\eqref{eq:rozvojela}, one can use results of~\cite{ADMP} in order to describe the distances between consecutive $(-\beta)$-integers and
substitution properties of the corresponding infinite word.

\begin{prop}\label{p:mezeryCP1}
Let $\beta>1$ be the root of $x^2-mx-n$, $m\geq n\geq 1$. Let $\Z_{-\beta}=\{t_j\mid j\in\Z\}$, where $t_0=0$ and $t_j<t_{j+1}$.
Then the distances between consecutive $(-\beta)$-integers take values $t_{j+1}-t_j \in \{\Delta_0^-,\Delta_1^-\}$, where
$$
\Delta_0^-=1 \quad \text{ and }\quad
\Delta_1^-=\begin{cases}
\frac{m}\beta & \text{ if } \ m = n\,,\\
1+\frac{n}\beta & \text{ if } \ m< n\,.\\
\end{cases}
$$
\end{prop}

\pfz
If $m>n$, then the string $d_{-\beta}(l)=m(m-n)^\omega$ satisfies conditions of Theorem~18 from~\cite{ADMP} and the statement follows. For $m=n$, we have $d_{-\beta}(l)=m0^\omega$ and we use Theorem~21 from~\cite{ADMP}.
\pfk

Similarly as in case of positive base, one can code the sequence of $(-\beta)$-integers by an infinite word, say $u_{-\beta}$, over the alphabet $\{0,1\}$, where the letter 0 is used if distance $\Delta_0^-$ occurs, and $1$ otherwise. Since negative numbers may have integer $(-\beta)$-expansion with non-negative digits, the infinite word is pointed bidirectional, where we denote the position of 0 by a delimiter $|$. Formally, we define
\begin{equation}\label{eq:defu}
u_{-\beta}=\cdots u_{-2}u_{-1}|u_{0}u_1u_2\cdots\in\{0,1\}^\Z\,,\qquad \text{where \ $u_i=X$ \ if \ $t_{i+1}-t_i=\Delta_X^-$.}
\end{equation}
For further use, let us denote the suffix of $u_{-\beta}$ starting at the delimiter by $u_{-\beta}^+=u_0u_1u_2\cdots$.

The results of~\cite{ADMP} can be used to derive the following statement about invariance of $u_{-\beta}$ under an antimorphism.

\begin{prop}\label{p:antimorfCP1}
Let $\beta>1$ be the root of $x^2-mx-n$, $m\geq n\geq 1$. Then the infinite word $u_{-\beta}$ is a fixed point of the antimorphism
$\overleftarrow{\varphi_{\text{-}\beta}}:\{0,1\}^*\to\{0,1\}^*$, given by
$$
\begin{aligned}
\overleftarrow{\varphi_{\text{-}\beta}}(0)&=0^{m-1}1,\\
\overleftarrow{\varphi_{\text{-}\beta}}(1)&=0^{m+n-1}1,
\end{aligned} \quad \text{ for $m>n$, } \quad\text{and}\quad
\begin{aligned}
\overleftarrow{\varphi_{\text{-}\beta}}(0)&=0^{m}1,\\
\overleftarrow{\varphi_{\text{-}\beta}}(1)&=0^m,
\end{aligned} \quad \text{ for $m=n$.}
$$
In particular, we have
$$
{u_{-\beta}}=\lim_{n\rightarrow +\infty}\overleftarrow{\varphi_{\text{-}\beta}}^n(1)|\overleftarrow{\varphi_{\text{-}\beta}}^n(0)\,.
$$
\end{prop}

For the class of quadratic Pisot numbers $\beta$, which are not covered by Proposition~\ref{p:mezeryCP1}, we have by~\eqref{eq:rozvojelb} that $d_{-\beta}(l) = \big((m-1)n\big)^\omega$. For such a case, results of~\cite{ADMP} cannot be used neither for deriving the distances, nor for determining invariance under an antimorphism, thus we provide our own demonstration.
For that, we need to consider the admissibility condition~\eqref{eq:ISadmiss} for this particular case, which reads as follows: A string of integer digits $y_1y_2y_3\cdots$
is a $(-\beta)$-expansion of a number $x\in[-\frac\beta{\beta+1},\frac1{\beta+1})$ if and only if every suffix $y_iy_{i+1}y_{i+2}\cdots$ satisfies
$$
\big((m-1)n\big)^\omega  \preceq_{\text{\tiny alt}} y_iy_{i+1}y_{i+2}\cdots \prec_{\text{\tiny alt}} 0\big((m-1)n\big)^\omega\,.
$$
In particular, for strings ending in $0^\omega$, this results in requiring that the infinite word $y_1y_2y_3\cdots$ does not contain a forbidden factor from the set $\{(m-1)A \mid 0\leq A< n\}$.

The following proposition is stated with more details than previous results about distances in $\Z_\beta$ and $\Z_{-\beta}$, in order to facilitate the proof of Proposition~\ref{p:antimorfCP2}.

\begin{prop}\label{p:mezeryCP2}
Let $\beta>1$ be the root of $x^2-mx+n$, $m-2\geq n\geq 1$. Let $\Z_{-\beta}=\{t_j\mid j\in\Z\}$, where $t_0=0$ and $t_j<t_{j+1}$.
Then the distances between consecutive $(-\beta)$-integers take values $t_{j+1}-t_j \in \{\Delta_0^-,\Delta_1^-\}$, where
$$
\Delta_0^-=1 \quad \text{ and }\quad
\Delta_1^-= 2-\frac{n}\beta\,.
$$
In particular, let $x<y$ be consecutive $(-\beta)$-integers.
\begin{enumerate}
    \item If $\langle x\rangle_{-\beta}=x_k\dots x_1 A\bullet\,,$ with $A\leq m-3\,,$ then $y-x=\Delta_0^-=1$.
    \item If $\langle x\rangle_{-\beta}=x_k\dots x_1(m-2)\bullet\,,$ then $y-x=\Delta_1^-=2-\frac n \beta$.
\end{enumerate}
\end{prop}

\pfz
 Since the natural order of reals corresponds to the alternate order of their $(-\beta)$-expansions, it suffices to find, for every $(-\beta)$-integer $x$
 with $\langle x\rangle_{-\beta}=x_k\dots x_1 A\bullet\,,$
 the $(-\beta)$-integer $y$ with the smallest alternately greater expansion. We do not consider strings ending with the digit $(m-1)$, because they are not admissible.
  \begin{enumerate}
  \item Let $A\leq m-3\,.$ Then $y=x+\Delta_0=x+1$, because
  $$
  \begin{array}{ccccccc}
  x&=&x_k&\dots& x_1&A&\bullet\\
  x+1&=&x_k&\dots& x_1&(A+1)&\bullet
  \end{array}
  $$
  \item Let $A=(m-2).$ Then $y=x+\Delta_1$, which can be justified as follows. The $(-\beta)$-expansion of $x$ is now of the form
  $$
  \langle x\rangle_{-\beta}=\dots X\,Y\,[(m-1)\,n]^k\,(m-2)\bullet
  $$
  where $k\in \mathbb N_0$ and $X\,Y\neq(m-1)\,n.$
    We distinguish two cases
    \begin{enumerate}
    \item Let $X\leq m-2,\,Y\geq 1$ or $X=m-1,\,Y\geq n+1$. Then
  $$
  \begin{array}{ccccccccc}
    x&=&\dots&X&Y&[(m-1)&n]^k&(m-2)&\bullet\\
    x+\Delta_1^-&=&\dots&X&(Y-1)&0&[(m-1)&n]^k&\bullet
  \end{array}
  $$
  \item Let $X\leq m-2,\,Y=0$. Then
  $$
  \begin{array}{ccccccccc}
    x&=&\dots&X&0&[(m-1)&n]^k&(m-2)&\bullet\\
    x+\Delta_1^-&=&\dots&X+1&m-1&[n&(m-1)]^k&n&\bullet
  \end{array}
  $$
  \end{enumerate}
  \end{enumerate}
\pfk

Since the distances between consecutive $(-\beta)$-integers take two values, denoted by $\Delta_0^-$ and $\Delta_1^-$, we may define
the pointed bidirectional infinite word $u_{-\beta}$ and one-directional infinite word $u_{-\beta}^+$ for the class of quadratic Pisot numbers $\beta$, roots of $x^2-mx+n$,
in the same way as before, namely by~\eqref{eq:defu}. In the following proposition we find the antimorphism fixing $u_{-\beta}$.

\begin{prop}\label{p:antimorfCP2}
Let $\beta>1$ be the root of $x^2-mx+n$, $m-2\geq n\geq 1$.
Then the infinite word $u_{-\beta}$ is a fixed point of the antimorphism
$\overleftarrow{\varphi_{\text{-}\beta}}:\{0,1\}^*\to\{0,1\}^*$, given by
\begin{align*}
\overleftarrow{\varphi_{\text{-}\beta}}(0)&=0^{m-2}1\\
\overleftarrow{\varphi_{\text{-}\beta}}(1)&=0^{m-2}10^{m-n-2}1
\end{align*}
In particular, we have
$$
{u_{-\beta}}=\lim_{n\rightarrow +\infty}\overleftarrow{\varphi_{\text{-}\beta}}^n(1)|\overleftarrow{\varphi_{\text{-}\beta}}^n(0)\,.
$$
\end{prop}

\pfz
From the definition, the set $\Zmb$ is self-similar, i.e.\  $(-\beta)\Zmb\subset\Zmb$. We will determine the antimorphism $\overleftarrow{\varphi_{\text{-}\beta}}$ under which $u_{-\beta}$ is invariant by showing that between consecutive points of $(-\beta)\Zmb$ (which lie also in $\Zmb$) we find always the same configuration of elements of $\Zmb$. More precisely, let $x<y$ be consecutive points in $\Zmb$. Then $-\beta y < -\beta x$ are consecutive points in $(-\beta)\Zmb$. We will show that there are only two types of the configurations $[-\beta x,-\beta y]\cap\Zmb$, according to whether $y-x=\Delta_0^-$ or $y-x=\Delta_1^-$. The sequence of distances in the configuration obtained for $y-x=\Delta_0^-$ determines the word $\overleftarrow{\varphi_{\text{-}\beta}}(0)$; similarly, the sequence of distances in the configuration obtained for $y-x=\Delta_1^-$ determines the word $\overleftarrow{\varphi_{\text{-}\beta}}(1)$. Recall that $\Delta_0^-=1$ is coded by the letter $0$ and $\Delta_1^-=2-\frac n\beta$ by the letter $1$.

%

Consider consecutive $(-\beta)$-integers $x,y$ such that $y=x+1=x+\Delta_0^-$. By Proposition~\ref{p:mezeryCP2}, the expansions of $x,y$ are of the form
$$
\begin{array}{ccccccc}
\langle x\rangle_{-\beta}&=&x_k&\dots& x_1&A&\bullet\\
\langle y\rangle_{-\beta}&=&x_k&\dots& x_1&(A+1)&\bullet
\end{array}
$$
Multiplying by $(-\beta)$ we have
$$
\begin{array}{cccccccc}
\langle -\beta x\rangle_{-\beta}&=&x_k&\dots& x_1&A&0&\bullet\\
\langle -\beta y\rangle_{-\beta}&=&x_k& \dots& x_1&(A+1)&0&\bullet
\end{array}
$$
and we have $-\beta x>-\beta y$. From Proposition~\ref{p:mezeryCP2}, we derive that $(m-2)$ right neighbours of the point $-\beta y$ are
$-\beta y+1$, $-\beta y+2$, \dots, $-\beta y+ m-2$.  There are no other elements of $\Zmb$ in the interval $(-\beta y, -\beta x)$.
We have $-\beta y+(m-2)\Delta_0^-+\Delta_1^-=-\beta x$, and thus we may put $\overleftarrow{\varphi_{\text{-}\beta}}(0)=0^{m-2}1$.


Let now $y=x+\Delta_1^-$. From the proof of Proposition~\ref{p:mezeryCP2}, we can derive that the number $-\beta y$ may have two possible expansions, namely
$$
\begin{array}{ccccccc}
\langle -\beta y\rangle_{-\beta} &=& \dots &(m-1)&n&0&\bullet\\
\langle -\beta y\rangle_{-\beta} &=& \dots &X&0&0&\bullet
\end{array}
$$
In both cases, the neighbouring $m-2$ elements of $\Zmb$ are again $-\beta y+1$, $-\beta y+2$, \dots, $-\beta y+ m-2$, and we have
$$
\begin{array}{ccccccc}
\langle -\beta y+(m-2)\Delta_0^-\rangle_{-\beta} &=& \dots &(m-1)&n&(m-2)&\bullet\\
\langle -\beta y+(m-2)\Delta_0^-\rangle_{-\beta} &=& \dots &X&0&(m-2)&\bullet
\end{array}
$$
where $X$ denotes an arbitrary digit in $\{0,1,\dots,m-2\}$. By Proposition~\ref{p:mezeryCP2}, the neighbour is
$$
\begin{array}{ccccccc}
\langle -\beta y+(m-2)\Delta_0^-+\Delta_1^-\rangle_{-\beta} &=& \dots &\dots&(m-1)&n&\bullet\\
\langle -\beta y+(m-2)\Delta_0^-+\Delta_1^-\rangle_{-\beta} &=& \dots &(X+1)&(m-1)&n&\bullet
\end{array}
$$
Necessarily, the $(m-n-2)$ right neighbours are obtained by adding $\Delta_0^-$.
The last one is of the form
$$
-\beta y + (m-2)\Delta_0^- + \Delta_1^- + (m-n-2)\Delta_0^- = -\beta x -\Delta_1^-\,.
$$
Hence there are no other points of $\Zmb$ in the interval $(-\beta y, -\beta x)$ and the prescription for
the antimorphism is $\overleftarrow{\varphi_{\text{-}\beta}}(1)=0^{m-2}10^{m-n-2}1$.
%
%
\pfk


As a consequence of the above descriptions of distances in $\Z_{-\beta}$ and antimorphisms fixing infinite words coding $\Z_{-\beta}$, we derive an analogue of Proposition~\ref{p:sturmianposit}, characterizing Sturmian words among $u_{-\beta}$, which is defined analogously to~\eqref{eq:defu} for any Ito-Sadahiro number $\beta$.

\begin{coro}\label{p:sturm}
Let $\beta>1$ be an Ito-Sadahiro number and let $u_{-\beta}$ be the bi-directional infinite word coding $\Z_{-\beta}$. Then $u_{-\beta}$ is Sturmian if and only if $\beta$ is a quadratic Pisot unit.
\end{coro}

\pfz
From Proposition~16 of~\cite{ADMP} it follows that the infinite word $u_{-\beta}$ is binary only if $\beta$ is quadratic. As mentioned above, the only quadratic Ito-Sadahiro numbers are quadratic Pisot numbers. For $\beta$ a quadratic Pisot unit, Proposition~\ref{p:negatCP1} identifies the set $\Z_{-\beta}$ as a cut-and-project set $\Sigma_\beta(\Omega)$ for certain interval $\Omega$. Propositions~\ref{p:mezeryCP1} and~\ref{p:mezeryCP2} state that the distances between consecutive elements of $\Z_{-\beta}$ take two values. By Theorem~\ref{t:zgumape}, the infinite word $u_{-\beta}$ is a bidirectional Sturmian word. On the other hand, if $\beta$ is a quadratic Pisot number, but not a unit, then we use the combinatorial characterization of Sturmian words as aperiodic balanced words to exclude that $u_{-\beta}$ is Sturmian.

Recall that words $u_{-\beta}$ are fixed points
$$
u_{-\beta}=\lim_{n\to\infty}\overleftarrow{\varphi_{\text{-}\beta}}^n(1)|\overleftarrow{\varphi_{\text{-}\beta}}^n(0)
$$
of antimorphisms given in Proposition~\ref{p:antimorfCP1} and~\ref{p:antimorfCP2}.
We distinguish three cases according to the prescription of the antimorphism.
\begin{enumerate}
\item Let $\beta^2=m\beta+m$, $m\geq 2$. Then $\overleftarrow{\varphi_{\text{-}\beta}}(0)=0^m1$, $\overleftarrow{\varphi_{\text{-}\beta}}(1)=0^m$. Necessarily, $u_{-\beta}$ contains the factors $10^m1$ and $0^{2m}$, which contradicts the balance property of Sturmian words.

\item Let $\beta^2=m\beta+n$, $2\leq n<m$. Then $\overleftarrow{\varphi_{\text{-}\beta}}(0)=0^{m-1}1$, $\overleftarrow{\varphi_{\text{-}\beta}}(1)=0^{m+n-1}1$. Then  $u_{-\beta}$ contains the factors $10^{m-1}1$ and $0^{m+n-1}$, which is a contradiction.

\item Let $\beta^2=m\beta-n$, $2\leq n \leq m-2$. Then $\overleftarrow{\varphi_{\text{-}\beta}}(0)=0^{m-2}1$, $\overleftarrow{\varphi_{\text{-}\beta}}(1)=0^{m-2}10^{m-n-2}1$. Then  $u_{-\beta}$ contains the factors $10^{m-n-2}1$ and $0^{m-2}$, which is again a contradiction.
\end{enumerate}
\pfk

\begin{pozn}
In fact, when showing that $u_{-\beta}$ is a Sturmian word for a quadratic Pisot unit $\beta$
one could avoid argumentation by cut-and-project scheme, providing thus an alternative proof of one implication in Corollary~\ref{p:sturm}. For that, if would be sufficient to demonstrate that $u_{-\beta}$ is a fixed point of a primitive Sturmian morphism.
A morphism $\sigma$ is called Sturmian, if for every Sturmian word $u$, the infinite word $\sigma(u)$ is Sturmian.
A morphism over a binary alphabet is Sturmian if and only if it is a composition of morphisms $E,\varphi,\widetilde\varphi$, where
\begin{center}
\begin{tabular}{rl}
\multirow{2}{*}{$E:$} & $0\mapsto 1$\\
 & $1\mapsto 0$
\end{tabular}
\begin{tabular}{rl}
\multirow{2}{*}{$\varphi:$} & $0\mapsto 01$\\
 & $1\mapsto 0$
\end{tabular}
\begin{tabular}{cl}
\multirow{2}{*}{$\widetilde\varphi:$} & $0\mapsto 10$\\
 & $1\mapsto 0$
\end{tabular}
\end{center}
The infinite word $u_{-\beta}$ is a fixed point of the morphism $\overleftarrow{\varphi_{\text{-}\beta}}^2$, where $\overleftarrow{\varphi_{\text{-}\beta}}$
is the antimorphism from Propositions~\ref{p:antimorfCP1} and~\ref{p:antimorfCP2}. We can justify that $\overleftarrow{\varphi_{\text{-}\beta}}^2$ is Sturmian  by providing its decomposition into morphisms $E,\varphi,\widetilde\varphi$. We have

\begin{enumerate}
\item Let $\beta^2=\beta+1$. Then $\beta=\tau=\frac12(1+\sqrt5)$, and $u_{-\beta}$ is a fixed point of the Sturmian morphism $\overleftarrow{\varphi_{\text{-}\beta}}^2=\varphi\circ\widetilde\varphi$, hence it is a Sturmian word.

\item Let $\beta^2=m\beta+1$, $1<m$. Then the second iteration $\overleftarrow{\varphi_{\text{-}\beta}}^2$ is a morphism with prescription
    $$
    0\mapsto 0^{m}1(0^{m-1}1)^{m-1}\,,\qquad 1\mapsto 0^{m}1(0^{m-1}1)^{m}\,.
    $$
    One can verify that
    $$
    \overleftarrow{\varphi_{\text{-}\beta}}^2 = (\varphi\circ E)^{m-1}\circ(E\circ\widetilde{\varphi})^m\circ\varphi\circ E\,.
    $$

\item Let $\beta^2=m\beta-1$. Then  $\overleftarrow{\varphi_{\text{-}\beta}}^2$ is given by the prescription
    $$
    0\mapsto 0^{m-2}10^{m-3}1(0^{m-2}1)^{m-2}\,,\qquad 1\mapsto 0^{m-2}10^{m-3}1(0^{m-2}1)^{m-2}0^{m-3}1(0^{m-2}1)^{m-2}\,.
    $$
    This morphism is Sturmian since
    $$
    \overleftarrow{\varphi_{\text{-}\beta}}^2 = (\varphi\circ E)^{m-3}\circ E\circ\widetilde{\varphi}\circ(\widetilde{\varphi}\circ E)^{m-2}\circ E\circ\widetilde{\varphi}\circ\varphi\circ E\,.
    $$
\end{enumerate}
\end{pozn}

\section{Relation of $\Z_\beta$ and $\Z_{-\beta}$}\label{sec:vztah}

One of the main purposes of this paper is demonstrating the relation of $\beta$- and $(-\beta)$-integers and thus continuing the study performed for the golden ratio $\tau$ in~\cite{MaVa}. There we have shown that
\begin{equation}\label{eq:vztahtau}
\Z_{-\tau}\cap \mathbb R^+=\Z_{\tau^2}^+\,.
\end{equation}
In fact, such an equality relating integers in base $-\beta$ and $\beta^2$ is exceptional, ensured by the fact that $\Z[\tau]=\Z[\tau^2]$. This allows that
\begin{equation}\label{eq:vztahtau2}
\Z_{\tau^2}\cap \mathbb R^+=(\Z_{\tau}+1)\cap \mathbb R^+\,.
\end{equation}
This can be seen for example from the fact that $\tau$ is a root of $x^2-x-1$, and so by the first statement of Proposition~\ref{p:positCP1},
$$
\Z_{\tau}^+ = \Sigma_\tau\big((-1,\tau)\big) \cap \mathbb R^+= \big\{x\in\Z[\tau] \,\big|\, x'\in (-1,\tau) \big\} \cap \mathbb R^+\,,
$$
and similarly, $\tau^2$ is a root of $x^2-3x+1$, and so by the second statement of Proposition~\ref{p:positCP1},
$$
\begin{aligned}
\Z_{\tau^2}^+ &=  \Sigma_{\tau^2}\big([0,\tau^2)\big) \cap \mathbb R^+ = \big\{x\in\Z[\tau^2]=\Z[\tau] \,\big|\, x'\in [0,\tau^2) \big\} \cap \mathbb R^+=\\
&=  \Sigma_{\tau}\big([0,\tau^2)\big) \cap \mathbb R^+=\Z_{-\tau}\cap \mathbb R^+\,.
\end{aligned}
$$
We obviously have $(-1,\tau)+1=(0,\tau^2)$, which according to statement 2 of Proposition~\ref{l:cap} and Proposition~\ref{p:mezeryCP2} justifies~\eqref{eq:vztahtau2} and~\eqref{eq:vztahtau}.

In~\cite{MaVa} we have used a combinatorial approach for proving~\eqref{eq:vztahtau}. Argumentation using the cut-and-project scheme as above can be extended to clarify the relation of $\Z_{-\beta}$ and $\Z_\beta$ for all quadratic Pisot units $\beta$ roots of $x^2-mx-1$, $m\geq 1$.

\begin{thm}\label{t:+-1}
Let $\beta>1$ be the root of $x^2-mx-1$, $m\geq 2$. Then
$$
\Z_\beta\cap\R^+=(\Z_{-\beta}\cup\beta\Z_{-\beta})\cap\R^+.
$$
\end{thm}

\pfz
By Proposition~\ref{p:positCP1} we have $\Z_\beta\cap\R^+=\Sigma(-1,\beta)\cap \R^+$. Since $\beta'=-\frac1\beta$, we have from Proposition~\ref{l:cap} that
$$
\Sigma_\beta(-1,\beta)=\Sigma_\beta(-1,0]\cup\Sigma_\beta[0,\beta)
=\Sigma_\beta\big(-\tfrac1\beta[0,\beta)\big)\cup\Sigma_\beta[0,\beta)=\beta\Sigma_\beta[0,\beta)\cup\Sigma_\beta[0,\beta).
$$
By Proposition~\ref{p:negatCP1} we have $\Z_{-\beta}=\Sigma[0,\beta)$, and the statement of the theorem follows.
\pfk

Note that $\Z_{-\beta}\cap\beta\Z_{-\beta}=\{0\}$, as can be seen from the proof of Theorem~\ref{t:+-1}.

If we try to obtain similar geometric relation for $\beta>1$ root of $x^2-mx+1$, $m\geq 3$, we do not succeed, for,
the closures of $\{x'\mid x\in\Z_\beta\}$, resp. $\{x'\mid x\in\Z_{-\beta}\}$ are shifted with respect to each other by a constant that does not belong to $\Z[\beta']=\Z[\beta]$. Even if we consider for example the set
$$
(\Z_\beta\setminus\beta\Z_\beta)\cap\R^+=\Sigma_\beta[1,\beta)\cap\R^+.
$$
which has the same interval length as the cut and project set
$$
\Z_{-\beta}=\Sigma_\beta\Big(\big(-\tfrac{\beta-1}{\beta+1},\beta\tfrac{\beta-1}{\beta+1}\big)\Big)\,,
$$
we cannot identify them as a shift one of the other. Nevertheless, every finite piece in $\Sigma_\beta\Big(\big(-\tfrac{\beta-1}{\beta+1},\beta\tfrac{\beta-1}{\beta+1}\big)\Big)$ can be found translated in $\Sigma_\beta[1,\beta)$
and vice versa. This suggests that we can still obtain an interesting combinatorial relation between $\Z_\beta$ and $\Z_{-\beta}$
when studying the corresponding infinite words $u_\beta$ and $u_{-\beta}$. In particular, we will study the language of these infinite words, thus comparing the finite pieces of $\Z_\beta$ and $\Z_{-\beta}$. In fact, such an approach allows us to give a relation between the $\beta$- and $(-\beta)$-integers not only for quadratic Pisot units $\beta$, but also for other quadratic Pisot numbers.

Since the language $\L(u_{-\beta})$ of $u_{-\beta}=\cdots u_{-2}u_{-1}|u_0u_1u_2\cdots$ is the same as the one of its suffix $u_{-\beta}^+=u_0u_1u_2\cdots$, it is natural to compare $\L(u_{\beta})$ and $\L(u_{-\beta}^+)$. The tool is to compare the morphisms under which the
words are invariant. For that, the notion of conjugation is useful. A morphism $\psi$ is a right conjugate of a morphism $\sigma$, if there is a finite word $w$, such that
$$
\sigma(a)w=w\psi(a)\quad\text{ for all }\ a\in\A^*\,.
$$
Incidence matrices $M_\sigma$ and $M_\psi$ of conjugate morphisms coincide. The opposite implication is obviously not satisfied.
We will use the following lemma, which appears to be a folklore, thus we include it without proof.

\begin{lem}\label{l:languagesconjugate}
Languages of fixed points of conjugate primitive morphisms coincide.
\end{lem}

In the terminology of dynamical systems, the fact that languages of two infinite words coincide means that they belong to the same subshift. However, it does not necessarily mean that one is a shift of the other.

\begin{thm}
Let $\beta>1$ be the root of $x^2-mx-m$, $m\geq 1$. Then $\L(u_\beta)=\L(u_{-\beta})$.
Moreover, if $m=1$, then $u_{-\beta}^+=0u_\beta$. When $m\geq 2$, then  neither $u_{-\beta}^+$ is a suffix of $u_\beta$,
nor $u_\beta$ is a suffix of $u_{-\beta}^+$.
\end{thm}

\pfz
By Proposition~\ref{p:antimorfCP1}, we have $\overleftarrow{\varphi_{\text{-}\beta}}(0)=0^m1$, $\overleftarrow{\varphi_{\text{-}\beta}}(1)=0^m$
and thus
$$
\overleftarrow{\varphi_{\text{-}\beta}}^2(0)=0^m(0^m1)^m\,,\qquad \overleftarrow{\varphi_{\text{-}\beta}}^2(1)=(0^m1)^m\,.
$$
The second iteration of the canonical substitution $\varphi_\beta$, $\varphi_\beta(0)=0^m1$, $\varphi_\beta(1)=0^m$ is given by
$$
\varphi^2(0)=(0^m1)^m0^m\,,\qquad \varphi^2(1)=(0^m1)^m\,.
$$
Therefore $\overleftarrow{\varphi_{\text{-}\beta}}^2$ is a right conjugate of $\varphi^2$ and the conjugation factor $w$ is equal to $w=(0^m1)^m$.

If $m=1$, then $\beta=\tau$ and the statement $u_{-\tau}^+=0u_\tau$ is found in~\cite{MaVa}.
On the other hand, if $m\geq 2$, we can exclude that $u_\beta$ is a suffix of $u_{-\beta}^+$ or vice versa.
Recall that letters $0,1$ code the same distances $\Delta_0^+=\Delta_0^-=1$, $\Delta_1^+=\Delta_1^-=\frac{m}{\beta}$ in $u_\beta$ and $u_{-\beta}^+$,  cf. Propositions~\ref{p:mezeryposit} and~\ref{p:mezeryCP1}. If $u_{-\beta}^+=wu_{\beta}$ or $wu_{-\beta}^+=u_{\beta}$ for some finite word $w\in\{0,1\}^*$, then
\begin{equation}\label{eq:c}
(\Z_{-\beta}\pm c)\cap\mathbb R^+=\Z_\beta^+\,,
\end{equation}
where $c=|w|_0\Delta_0+|w|_1\Delta_1\in\Z+\Z\frac{m}\beta=\Z[\beta]$. Taking the Galois images, one would have
$$
\{x'\mid x\in\Z_{-\beta}\cap\R^+\}\pm c' = \{x'\mid x\in\Z_\beta^+\}\,.
$$
However, this is not possible. We will show that by finding $\sup\{x'\mid x\in\Z_{-\beta}\cap\R^+\}$ and $\sup\{x'\mid x\in\Z_\beta^+\}$.
It is easy to realize that supremum $\sup\{x'\mid x\in\Z_\beta^+\}$ is approached by $\beta$-integers $x=\sum_{i=0}^ka_i\beta^i\in\Z[\beta]$ with $\beta$-expansion
$\langle x\rangle_\beta = m0m0\cdots m0m \bullet$. We have
$$
x'=\sum_{i=0}^ka_i{\beta'}^i = \sum_{i=0}^ka_i{\Big(-\frac{m}\beta\Big)}^i < m \sum_{i=0}^\infty {\Big(\frac{m}\beta\Big)}^{2i} = \frac{\beta^3}{\beta+m}\,.
$$
Similarly, one can compute $\sup\{x'\mid x\in\Z_{-\beta}\cap\R^+\}$, realizing that it is obtained considering $(-\beta)$-integers of the form
$\langle x\rangle_{-\beta} = mmm\cdots m \bullet$,
$$
x'=\sum_{i=0}^ka_i(-\beta')^i < m \sum_{i=0}^\infty {\Big(\frac{m}\beta\Big)}^{i} = \frac{m\beta}{\beta-m} = \beta^2\,.
$$
Thus the difference $\sup\{x'\mid x\in\Z_\beta^+\}-\sup\{x'\mid x\in\Z_{-\beta}\cap\R^+\}=\frac{\beta^3}{\beta+m}-\beta^2\notin\Z[\beta]$, and hence equality~\eqref{eq:c}
cannot be valid for any $c\in\Z[\beta]$.
\pfk

One cannot expect that morphisms corresponding to $u_\beta$ and $u_{-\beta}$ will be conjugated for  quadratic Pisot $\beta$ other than roots of $x^2-m x -m$. This is because the incidence matrices of morphisms $\varphi_\beta^2$ and $\overleftarrow{\varphi_{\text{-}\beta}}^2$ do not coincide, compare Proposition~\ref{p:substposit}
with Propositions~\ref{p:antimorfCP1} and~\ref{p:antimorfCP2}. This is also in agreement with the fact that the distances in $\Z_\beta$ and $\Z_{-\beta}$ coincide only
if $\beta$ is a root of $x^2-m x -m$. For other quadratic Pisot numbers $\beta$, we have interestingly
$$
\begin{aligned}
\Delta_0^+ &=\Delta_0^-=1,\quad \Delta_1^+=\frac{n}{\beta} = \Delta_1^--1\quad &\ \text{ if }\ \beta^2=mx+n,\ m>n\geq 1\,,\ \ \quad\\
\Delta_0^+ &=\Delta_0^-=1,\quad \Delta_1^+=1-\frac{n}{\beta} = \Delta_1^--1\quad&\text{if }\ \beta^2=mx-n,\ m-2\geq n\geq 1\,,
\end{aligned}
$$
which means that $\Delta_1^-=\Delta_0^++\Delta_1^+$. It is therefore appealing to see what happens when adding to $\Z_{-\beta}$ certain points which would split every distance $\Delta_1^-$ into two distances of lengths $\Delta_0^+=\Delta_0^-=1$ and $\Delta_1^+$. Theorem~\ref{t:+-1} shows that at least for one class of considered numbers, this splitting on the positive half-line gives the same set as $\Z_\beta^+$. Note, however, that we can choose the order in which we split the long distance into two.

Realize that performing this procedure on the geometric representation of $\Z_{-\beta}$ corresponds to applying a suitable morphism on the infinite word $u_{-\beta}$.
Define $\pi, \tilde{\pi}:\{0,1\}^*\to\{0,1\}^*$ by
\begin{equation}\label{eq:pi}
\text{\begin{tabular}{rl}
\multirow{2}{*}{$\pi:$} & $0\mapsto 0$\\
 & $1\mapsto 10$
\end{tabular}
\qquad{ and }\qquad
\begin{tabular}{rl}
\multirow{2}{*}{$\tilde{\pi}:$} & $0\mapsto 0$\\
 & $1\mapsto 01$.
\end{tabular}
}
\end{equation}
As we shall prove, application of $\pi$, resp. $\tilde{\pi}$ on the infinite word $u_{-\beta}$ leads to an infinite word with the same language as $u_{\beta}$.
For showing that, we will use the following easy statement.

\begin{lem}\label{l:pevbod}
Let $\sigma, \psi, \pi$ be morphisms such that $\sigma:\A^*\to\A^*$,  $\psi:\B^*\to\B^*$, $\pi: \A^*\to\B^*$
and $\pi\circ\sigma = \psi\circ\pi$. If $u$ is a fixed point of $\sigma$, then $v:=\pi(u)$ is a fixed point of $\psi$.
\end{lem}

\pfz
We have $v=\pi(u)=\pi\big(\sigma(u)\big)=\psi\big(\pi(u)\big)=\psi(v)$.
\pfk

\begin{thm}\label{t:konjugrozsekane}
Let $\pi$ and $\tilde{\pi}$ be morphisms as in~\eqref{eq:pi}.

Let $\beta>1$ be the root of $x^2-mx-n$, $m> n\geq 1$.
Then,
$$
\L(u_\beta) =  \L\big(\tilde{\pi}(u_{-\beta})\big)\,.
$$
Moreover, $u_\beta=\pi(u_{-\beta}^+)$ for $n=1$.

Let $\beta>1$ be the root of $x^2-mx+n$, $m-2\geq n\geq 1$.
$$
\L(u_\beta) =  \L\big({\pi}(u_{-\beta})\big)\,.
$$
\end{thm}

\pfz
Let $\beta$ be the quadratic Pisot number satisfying $\beta^2=m\beta+n$, $m> n\geq 1$. Define a morphism
$\psi:\{0,1\}^*\to\{0,1\}$ by
$$
\begin{tabular}{rl}
\multirow{2}{*}{$\psi:$} & $0\mapsto 0^{m+n}1(0^m1)^{m-1}$\\
 & $1\mapsto (0^m1)^n$
\end{tabular}
$$
One can easily verify that the morphism $\psi$ satisfies $\tilde{\pi}\circ\overleftarrow{\varphi_{\text{-}\beta}}^2 = \psi\circ\tilde{\pi}$.
Therefore by Lemma~\ref{l:pevbod}, $u_{-\beta}$ is a fixed point of $\psi$. We can also verify that $\psi$ is a right conjugate of
$\varphi_\beta^2$,  which is by Proposition~\ref{p:substposit} equal to
$$
\begin{tabular}{rl}
\multirow{2}{*}{$\varphi_\beta^2:$} & $0\mapsto (0^m1)^m0^n$\\
 & $1\mapsto (0^m1)^n$
\end{tabular}
$$
The conjugacy factor $w$ is equal to $w=(0^m1)^m $. By Lemma~\ref{l:languagesconjugate}, the infinite words $u_\beta^+$ and $u_{-\beta}$
have the same language.

Similarly, for the quadratic Pisot number satisfying $\beta^2=m\beta-n$, $m-2\geq n\geq 1$, define
$\psi:\{0,1\}^*\to\{0,1\}$ by
$$
\begin{tabular}{rl}
\multirow{2}{*}{$\psi:$} & $0\mapsto (0^{m-2}1)0^{m-n-1}1(0^{m-1}1)^{m-2}0$\\
 & $1\mapsto (0^{m-2}1)0^{m-n-1}1(0^{m-1}1)^{m-n-2}0$
\end{tabular}
$$
The morphism $\psi$ satisfies ${\pi}\circ\overleftarrow{\varphi_{\text{-}\beta}}^2 = \psi\circ{\pi}$ and therefore $u_{-\beta}$ is a fixed point of $\psi$.
We verify that $\varphi_\beta^2$ is a right conjugate of
$\psi$,  which is by Proposition~\ref{p:substposit} equal to
$$
\begin{tabular}{rl}
\multirow{2}{*}{$\varphi_\beta^2:$} & $0\mapsto (0^{m-1}1)^{m-1}0^{m-n-1}1$\\
 & $1\mapsto (0^{m-1}1)^{m-n-1}0^{m-n-1}1$
\end{tabular}
$$
The conjugacy factor $w$ is equal to $w=0^{m-2}10^{m-n-1}1 $. Thus by Lemma~\ref{l:languagesconjugate}, we derive $\L(u_\beta) =  \L\big({\pi}(u_{-\beta})\big)$.

It remains to justify that for $\beta^2=m\beta+1$ we have equality $u_\beta=\pi(u_{-\beta}^+)$.
In this case
\begin{equation}
\text{\begin{tabular}{rl}
\multirow{2}{*}{$\varphi_\beta^2:$} & $0\mapsto (0^m1)^m0$\\
 & $1\mapsto 0^m1$
\end{tabular}
\qquad{ and }\qquad
\begin{tabular}{rl}
\multirow{2}{*}{$\overleftarrow{\varphi_{\text{-}\beta}}^2:$} & $0\mapsto 0^m1(0^{m-1}1)^{m-1}$\\
 & $1\mapsto 0^m1(0^m1)^m$
\end{tabular}
}
\end{equation}
and therefore one can verify that ${\pi}\circ\overleftarrow{\varphi_{\text{-}\beta}}^2 = \varphi_\beta^2\circ{\pi}$. By Lemma~\ref{l:pevbod},
$\pi(u_{-\beta}^+)$ is a fixed point of $\varphi_\beta^2$, but $\varphi_\beta^2$ has only one fixed point, namely $u_\beta=\lim_{j\to\infty} \varphi_\beta^j(0)$.
\pfk

\section{Addition of $(-\beta)$-integers}\label{sec:group}

In the following, let us focus on the arithmetical aspect of $(-\beta)$-integers. For non-integer base $-\beta$, the set $\Z_{-\beta}$ is not closed under addition and subtraction. One can nevertheless define the operation $\oplus$ of addition on $\Z_{-\beta}$ as
$t_j\oplus t_k = t_{j+k}$, with the neutral element $t_0=0$ and opposite element $\ominus t_j= t_{-j}$,
so that $\Z_{-\beta}$ be an additive group isomorphic to $\Z$. Such a definition is of course possible for any countable set.
Here we show that for every quadratic Pisot number $\beta$, the `sum' $t_j\oplus t_k$ yields always a result which is not far from the ordinary sum of real numbers $t_j$ and $t_k$. More formally, we show that $t_j\oplus t_k - (t_j+t_k)$ is bounded independently on $j$ and $k$, and we determine the possible values of the distance for quadratic Pisot units. We also show that the operation $\oplus$ is compatible with ordinary addition of real numbers, namely that whenever the result
$t_j+t_k$ is a $(-\beta)$-integer, then $t_j\oplus t_k =t_j+t_k$.
Note that for the case of positive base $\beta$, where $\beta$ is a quadratic Pisot unit, this was done in~\cite{ElFrGaVG} by a technical arithmetic study. Here we choose a combinatorial approach which allows us to derive a general result about compatibility.
%
%
%
We will show that for any discrete set $\Sigma=\{t_j\mid j\in\Z\}\subset\R$, $t_j<t_{j+1}$, with finitely many distances between consecutive points,
i.e. values $t_{j+1}-t_j$, one can define a binary operation $\oplus$ so that $(\Sigma,\oplus)\simeq(\Z,+)$, which is compatible with addition in $\R$
if the values of the distances are linearly independent over $\Q$.

\begin{thm}\label{t:compat}
Let $u=\cdots u_{-2}u_{-1}|u_0u_1u_2\cdots$ be a bidirectional infinite word over a finite alphabet $\A=\{0,1,\dots,d-1\}$. Let
$\Delta_0$, \dots, $\Delta_{d-1}$ be positive real values linearly independent over $\Q$.
Set $t_0=0$. For $k\geq 1$, define $t_k=\sum_{i=0}^{d-1}|w|_i\Delta_i$, where $w=u_0u_1\cdots u_{k-1}$. For $k\leq -1$, define
$t_k=-\sum_{i=0}^{d-1}|w|_i\Delta_i$, where $w=u_{-k}u_{-k+1}\cdots u_{-1}$. On the set $\Sigma=\{t_k\mid k\in\Z\}$ define a binary operation $\oplus$,
by $t_j\oplus t_k = t_{j+k}$, for $j,k\in\Z$. Then $(\Sigma,\oplus)$ is isomorphic to $(\Z,+)$ and $\oplus$ is compatible with addition in $\R$.
\end{thm}

\pfz
We have to show that whenever $t_j+t_k=t_l$, then $l=j+k$. For the contradiction, let us assume that $t_j+t_k=t_l,$ where $l\neq j+k.$ Suppose first that $0\leq j\leq k$ which obviously implies $l\geq 0.$ Let $t_j$ correspond to a prefix $w,\ t_k$ to a prefix $v$ and $t_l$ to a prefix $z$ of the word $u_0u_1\cdots$. We have
\begin{equation}\label{rovnice}t_j+t_k-t_l=0 \qquad\Longleftrightarrow\qquad \sum_{i=0}^{d-1}(|w|_i+|v|_i-|z|_i)\Delta_i=\sum_{i=0}^{d-1}(|wv|_i-|z|_i)\Delta_i=0.\end{equation}
The last equality of~\eqref{rovnice} is a nontrivial combination of $\Delta_i$'s which follows from $|wv|=j+k\neq l=|z|.$
Here we have a contradiction with the linear independence of  $\Delta_0,\dots,\Delta_{d-1}$ over $\mathbb Q.$

The same argument can be applied when considering indices $j\leq 0\leq k$ or $j\leq k\leq 0$.
\pfk

\begin{figure}[ht]
{
\begin{center}
\setlength{\unitlength}{3pt}
\linethickness{0.7pt}
\begin{picture}(130,20)
\put(5,10){\line(1,0){120}}
\put(10,9){\line(0,1){2}}
\put(50,9){\line(0,1){2}}
\put(70,9){\line(0,1){2}}
\put(115,9){\line(0,1){2}}
\put(110,9){\line(0,1){2}}
%
\put(10,13){$\overbrace{\hspace*{120pt}}^\text{\normalsize $w$}$}
\put(10,3.5){$\underbrace{\hspace*{180pt}}_\text{\normalsize $v$}$}
\put(70,13){$\overbrace{\hspace*{135pt}}^\text{\normalsize $w'$}$}
\put(6,5){$t_0=0$}
\put(49,5){$t_j$}
\put(69,5){$t_k$}
\put(115,5){$t_{j+k}$}
\put(101,5){$t_{j}+t_{k}$}
\end{picture}
\end{center}
}
\caption{Addition $t_j+t_k$ in $\Sigma$.}
\label{f:sigma}
\end{figure}
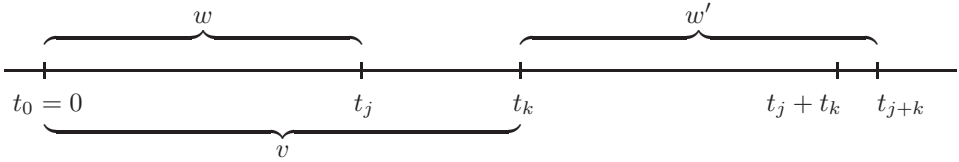

Let us derive the possible outcomes of $t_j+t_k-t_{j+k}.$
For $0\leq j\leq k$, let $t_j$ correspond to the prefix $w=u_0\dots u_{j-1}$, let $t_k$ correspond to the prefix $v=u_0\dots u_{k-1}$ and let $t_{j+k}$ correspond to the prefix $z=vw'$, where $w'=u_k\dots u_{j+k-1}$, see Figure~\ref{f:sigma}.
Then we have
\begin{equation}\label{eq:rozdil1}
  t_j+t_k-t_{j+k}=\sum_{i=0}^{d-1}(|w|_i+|v|_i-|vw'|_i)\Delta_i=\sum_{i=0}^{d-1}(|w|_i-|w'|_i)\Delta_i.
\end{equation}
The same looking formula can be derived for other cases $j\leq k\leq 0$ or $j\leq 0\leq k,$ only the factors
$w$ and $w'$ are dependent on the case, however, they are always of the same length.

\bigskip
Equality~\eqref{eq:rozdil1} states that the distance of $t_j+t_k$ from $t_j\oplus t_k=t_{j+k}$ depends on the difference of the number of occurrences of a specific letter in two factors $w,w'$ of the same length. This is captured by the notion of generalized balance defined in~\cite{adamczewski}. An infinite word $u$ over an alphabet $\A$ is called $C$-balanced for a $C\in\N$, if
$$
\big||w|_a-|w'|_a\big|\leq C\,,\quad\text{for every } a\in\A \text{ and every pair } w,w'\in\L(u),\ |w|=|w'|\,.
$$
We sometimes say that $u$ has bounded balances, if there exists $C<+\infty$, such that $u$ is $C$-balanced.
Obviously, $u$ is balanced if it is $1$-balanced.

Using~\eqref{eq:rozdil1}, we derive that if the set $\Sigma$ is defined by an infinite word with bounded balances, then the result of $t_j\oplus t_k$
is not far from the actual sum of $t_j$ and $t_k$.

\begin{coro}\label{c:balance}
Let $u=\cdots u_{-2}u_{-1}|u_0u_1u_2\cdots$ be a bidirectional infinite word over a finite alphabet $\A=\{0,1,\dots,d-1\}$ with bounded balances.
Define $(\Sigma,\oplus)$ as in Theorem~\ref{t:compat}. Then $t_j+t_k- (t_{j}\oplus t_k)$ is bounded independently of $j,k\in\Z$.
\end{coro}

Note that when $\Sigma$ is coded by a binary word over the alphabet $\A=\{0,1\},$
 we can write $|w|_1=|w|-|w|_0$. In this case,~\eqref{eq:rozdil1} can be rewritten as
\begin{equation}\label{eq:rozdilbinary}
t_j+t_k- t_{j+k}= \big(|w|_0-|w'|_0\big)\Delta_0 +\big(|w|_1-|w'|_1\big)\Delta_1 = \big(|w|_0-|w'|_0\big)(\Delta_0-\Delta_1)\,.
\end{equation}
Similarly, we would get
\begin{equation}\label{eq:rozdilbinary2}
- t_{-j}-t_j= \big(|w|_0-|w'|_0\big)(\Delta_0-\Delta_1) \,,
\end{equation}
where $w=u_{-j}u_{-j+1}\cdots u_{-1}$, $w'=u_0u_1\cdots u_{j-1}$  for $j\geq 0$ and vice versa otherwise.

\bigskip
Let us apply Theorem~\ref{t:compat} and Corollary~\ref{c:balance} to $\beta$- and $(-\beta)$-integers in case of quadratic Pisot $\beta$.
For a positive base, we know that the infinite word $u_\beta$ has bounded balances, they have been determined explicitly
in~\cite{balanceNonsimple,balanceSimple}. It is not difficult to prove that also the infinite word $u_{-\beta}$ has bounded balances, if we use the
result of Adamczewski~\cite{adamczewski} about balance properties of fixed points of morphisms. He shows that boundedness/unboundedness of balances is
decided by the second (in modulus) eigenvalue of the incidence matrix of the morphism.

\begin{lem}\label{l:balance}
For every quadratic Pisot $\beta$ there exist a constant $C$ such that the infinite word $u_{-\beta}$ is $C$-balanced.
If $\beta$ is a unit, then $C=1$. 
\end{lem}

\pfz
For any quadratic Pisot $\beta$, the infinite word $u_{-\beta}$ is a fixed point of a substitution, which is
a second iteration of the antimorphism $\overleftarrow{\varphi}$ given in Propositions~\ref{p:antimorfCP1} and~\ref{p:antimorfCP2}.
The incidence matrix of the substitution has $\beta^2$ as its dominant eigenvalue. Necessarily, the second eigenvalue is $\beta'^2$, which is in modulus smaller than 1.
Consequently, by~\cite{adamczewski}, the fixed point of $\overleftarrow{\varphi}_\beta^2$ has bounded balances.
If $\beta$ is a unit, then Proposition~\ref{p:sturm} states that the infinite word $u_{-\beta}$ is a Sturmian word, and thus it is balanced, i.e. $C=1$.
\pfk

Relations \eqref{eq:rozdilbinary} and~\eqref{eq:rozdilbinary2} describe the possible outcomes of $t_j\oplus t_k - (t_j+t_k)$ and $\ominus t_j - (-t_{j})$.
By Lemma~\ref{l:balance}, there exists a constant $C$ such that for any pair of factors $w,v\in{\mathcal L}(u_{-\beta})$
we have $\big||v|_0-|w|_0\big|\leq C$. This, together with~\eqref{eq:rozdilbinary} and~\eqref{eq:rozdilbinary2} shows that
$t_j\oplus t_k - (t_j+t_k)$ and $\ominus t_j - (-t_{j})$ are bounded independently of $j,k$, and that when
$\beta$ is a unit, then $C=1$ and thus
\begin{equation}\label{eq:rozdilunit}
t_j\oplus t_k - (t_j+t_k), \ominus t_j - (-t_{j}) \in\big\{0,\pm (\Delta_0-\Delta_1) \big\}\,.
\end{equation}

\begin{coro}
Let $\beta>1$ be the root of $x^2-mx+1$, $m\geq 3$.
Then
$$
t_j\pm t_k = t_{j\pm k} + \{-\eta,0,\eta\},\qquad\text{where $\eta=1-\frac1{\beta}$}.
$$
\end{coro}

For the second class of quadratic Pisot units, one can observe even stronger properties. This is given by the fact that the prefixes of $u_{-\beta}^+=u_0u_{1}u_3\dots$ are `light', or in other words, contain at least as many letters 0 as any factor of $u_{-\beta}$ of the same length. In fact, by Corollary~\ref{p:sturm}, $u_{-\beta}$ is a Sturmian word, i.e.\ it is balanced. Hence for any given length $k\in\N$, the number $|w|_0$ in factors $w\in\L(u_{-\beta})$, $|w|=k$, takes only two possible values. Factors having the larger number of $0$'s are called light, the other ones are called heavy.

\begin{lem}
Let $\beta>1$ be the root of $x^2-mx-1$, $m\geq 1$. Let $u_{-\beta}=\cdots u_{-2}u_{-1}|u_0u_1u_2\cdots $ be the infinite word coding the set of $(-\beta)$-integers.
Let $w=u_0u_1\cdots u_{n-1}$ and let $w'$ be any factor in $\L(u_{-\beta})$ of length $|w'|=|w|=n$.
Then $|w|_0-|w'|_0\in\{0,1\}$. In other words, every prefix of $u_{-\beta}^+$ is light.
\end{lem}

\pfz
When $\beta>1$ is the root of $x^2-mx-1$, $m\geq 1$, we use Proposition~\ref{p:negatCP1}, which states that $\Z_{-\beta}=\Sigma_\beta\big([0,\beta^2)\big)$ for $m=1$
and $\Z_{-\beta}=\Sigma_\beta\big([0,\beta)\big)$ otherwise. Consider first the case $m>1$.
The infinite word $u_{-\beta}^+$ codes the non-negative part of $\Z_{-\beta}$, and in fact, coincides with the Sturmian word arising as coding of 0 in the exchange of intervals $[0,\beta-1)$, $[\beta-1,1)$. By rescaling, it can be viewed also as a coding of 0 under a rotation by the angle $\alpha=\frac1\beta$ on the unit interval, thus
$u_{-\beta}^+$ is in fact what is called a Sturmian word with null intercept. As a consequence of the proof of Theorem~8 in~\cite{RiSaVa}, any prefix of such a Sturmian word is light. \pfk

\begin{prop}\label{p:soucet+}
Let $\beta>1$ be the root of $x^2-mx-1$, $m\geq 1$. Then
$$
t_j+t_k = t_{j+k} + \{0,\xi\}\,,\qquad t_j-t_k = t_{j-k} - \{0,{\xi}\}\qquad\text{for all } j,k\in\Z
$$
and
$$
-t_j= t_{-j} - \xi\quad \text{for } j\in\Z\setminus\{0\},
$$
where
$$
\xi=\Delta_0^--\Delta_1^-=\begin{cases}
-\frac1{\beta}  &\text{for }m\geq 2\,,\\
\frac1{\beta^2}&\text{for }m= 1\,.
\end{cases}
$$
Moreover, if $m\geq 2$, then $t_j\oplus t_k$ is the closest $(-\beta)$-integer to $t_j+t_k.$
\end{prop}

\pfz
Let us first prove $-t_j=t_{-j}-\xi$ for all $j\in\Z\setminus\{0\}.$
According to \eqref{eq:rozdilunit} we have that $-t_j=t_{-j}+\{0,\pm(\Delta_0^--\Delta_1^-)\},$ where $\Delta_0^--\Delta_1^-=\xi.$
We take the Galois images of those three possible cases, namely
\begin{equation}\label{eq:carky}
-t'_j=t'_{-j},\qquad -t'_j=t'_{-j}+\xi'\qquad\text{or}\qquad -t'_j=t'_{-j}-\xi'\,,
\end{equation}
where $\xi'=\beta^2$ and $\xi'=\beta$ for $m=1$ and $m\geq2$ respectively.
Proposition~\ref{p:negatCP1} states that $t_{j}\in\Zmb$ implies $t'_{j}\in\Omega$, where $\Omega=[0,\beta^2)$ for $m=1$ and $\Omega=[0,\beta)$ otherwise (note that $0$ corresponds to $t_0=0$ which is not considered in this statement).
Therefore $-t'_j\in-\Omega$ and substituting together with $t'_{-j}\in\Omega$ into \eqref{eq:carky} we get that the only possible option is $-t'_j=t'_{-j}-\xi'.$

\bigskip
Let $0\leq j\leq k.$ Using \eqref{eq:rozdil1} and the fact that the factor $w=u_0\dots u_{j-1}$ is light, we obtain $t_j+t_k- t_{j+k}=(|w|_0-|w'|_0)(\Delta_0^--\Delta_1^-)\in\{0,\xi\}.$ When $k\leq j\leq 0$ we can write
\begin{equation}\label{eq:ominusovani}
-(t_j+t_k-t_{j+k})=t_{-j}-\xi+t_{-k}-\xi-t_{-j-k}+\xi=
\underbrace{t_{-j}+t_{-k}-t_{-j-k}}_{\in\{0,\xi\}}-\xi\in\{-\xi,0\}.
\end{equation}
The first equality follows from $-t_j=t_{-j}+\xi$ and then we used already proven case $0\leq j\leq k.$
Thus is holds $t_j+t_k-t_{j+k}\in\{0,\xi\}.$

\begin{figure}[ht]\label{f:odcitani}
{
\begin{center}
\setlength{\unitlength}{3pt}
\linethickness{0.7pt}
\begin{picture}(130,20)
\put(5,10){\line(1,0){120}}
\put(10,9){\line(0,1){2}}
\put(50,9){\line(0,1){2}}
\put(70,9){\line(0,1){2}}
\put(115,9){\line(0,1){2}}
\put(55,9){\line(0,1){2}}
%
\put(10,13){$\overbrace{\hspace*{120pt}}^\text{\normalsize $w'$}$}
\put(50,3.5){$\underbrace{\hspace*{60pt}}_\text{\normalsize $z$}$}
\put(70,13){$\overbrace{\hspace*{135pt}}^\text{\normalsize $w$}$}
\put(9,5){$t_j$}
\put(48,5){$t_{j+k}$}
\put(53,13){$t_{j}+t_{k}$}
\put(66,5){$t_0=0$}
\put(115,5){$t_{k}$}
\end{picture}
\end{center}
}
\caption{Addition $t_j+t_k$ in the case $j\leq 0\leq k$ with $|j|\geq k$.}
\end{figure}
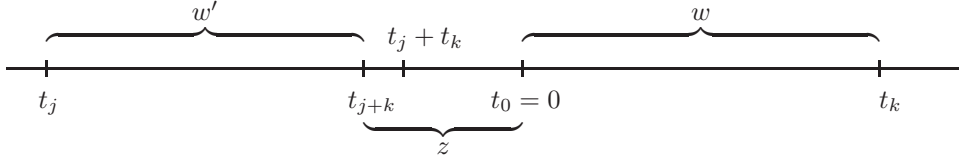

\bigskip
 Consider the case $j\leq 0\leq k$ with $|j|\geq k$. Let $w=u_{0}u_1\dots u_{k-1}$, $w'=u_{j}u_{-j+1}\dots u_{j+k-1}$, and $z=u_{j+k}u_{j+k+1}\dots u_{-1}$, i.e.\
 $t_k$ correspond to the prefix $w$, $t_j$ correspond to the suffix $v=w'z$, and $t_{j+k}$ to the suffix $z$, see Figure~\ref{f:odcitani}. Then we have
$$
t_j+t_k-t_{j+k}=\sum_{i=0}^1(-|w'z|_i+|w|_i+|z|_i)\Delta_i^-=(|w|_0-|w'|_0)(\Delta_0^--\Delta_1^-)=(|w|_0-|w'|_0)\xi.
$$
Since $w$ is light and factors $w,w'$ are of the same length, we have $|w|_0-|w'|_0\in\{0,1\}.$

\bigskip
When $j\leq 0\leq k$ with $|j|<k$ we use the same approach as in \eqref{eq:ominusovani} to get the statement.

\bigskip
Substituting $j=p,\ k=q-p,$ into $t_j+t_k=t_{j+k}+\{0,\eta\}$, one gets $t_p+t_{q-p}=t_q+\{0,\eta\}$, and thus $t_q-t_p=t_{q-p}-\{0,\eta\}$.

\bigskip

If $\beta>1$ is the root of $x^2-mx-1$, $m\geq 2$, then the distances between consecutive $(-\beta)$-integers take values
$\Delta_0^-=1$ and $\Delta_1^-=1+\frac1\beta$, and thus $\Delta_0^--\Delta_1^-=-\frac1\beta$. Moreover, since $\lfloor\beta\rfloor=m\geq 2$, we have
$\frac1\beta<\frac12$. Therefore both $t_j\oplus t_k - (t_j+t_k)$ and  $\ominus t_j - (-t_{j})$ are in modulus smaller
than $\frac12<\min(\Delta_0^-/2,\Delta_1^-/2)$.
\pfk

In the above proposition, we show that $t_j\oplus t_k$ is the closest $(-\beta)$-integer to $t_{j+k}$ only for $\beta$ root of $x^2-mx-1$ with $m\geq 2$. We can show that for $m=1$, i.e.\ when $\beta$ is the golden ratio, then similar statement does not hold.

\begin{ex}
Let $\beta$ be the golden ratio, root of $x^2-x-1$. Then $t_j\oplus t_k$ is not always the closest $(-\beta)$-integer to $t_j+t_k$.
The situation is illustrated in Figure~\ref{f:tau}.

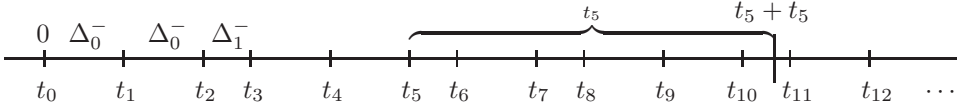
\begin{figure}[ht]
{
\begin{center}
\setlength{\unitlength}{3pt}
\linethickness{0.7pt}
\begin{picture}(130,15)
\put(5,5){\line(1,0){120}}
\put(10,4){\line(0,1){2}}
\put(20,4){\line(0,1){2}}
\put(30,4){\line(0,1){2}}
\put(36,4){\line(0,1){2}}
\put(46,4){\line(0,1){2}}
\put(56,4){\line(0,1){2}}
\put(62,4){\line(0,1){2}}
\put(72,4){\line(0,1){2}}
\put(78,4){\line(0,1){2}}
\put(88,4){\line(0,1){2}}
\put(98,4){\line(0,1){2}}
\put(104,4){\line(0,1){2}}
\put(114,4){\line(0,1){2}}
\put(102,2){\line(0,1){6}}
\put(97,10){$t_5+t_5$}
\put(56,7){$\overbrace{\hspace*{138pt}}^{t_5}$}
\put(13,7){$\Delta_0^-$}
\put(23,7){$\Delta_0^-$}
\put(31,7){$\Delta_1^-$}
\put(9,7){$0$}
\put(9,0){$t_0$}
\put(19,0){$t_1$}
\put(29,0){$t_2$}
\put(35,0){$t_3$}
\put(45,0){$t_4$}
\put(55,0){$t_5$}
\put(61,0){$t_6$}
\put(71,0){$t_7$}
\put(77,0){$t_8$}
\put(87,0){$t_9$}
\put(96,0){$t_{10}$}
\put(103,0){$t_{11}$}
\put(113,0){$t_{12}$}
\put(121,0){$\cdots$}
\end{picture}
\end{center}
}
\caption{Addition in $(-\beta)$-integers for $\beta=\frac12(1+\sqrt5)$.}
\label{f:tau}
\end{figure}

Since the distances between consecutive $(-\beta)$-integers are $\Delta_0^-=0$ and $\Delta_1^-=\frac1\beta$,
we can see in Figure~\ref{f:tau} that $t_5+t_5=2(4+\frac1\beta)=8+\frac2\beta$, $t_{10} = 7+\frac3\beta$ and $t_{11}=7+\frac4\beta$, so that
$$
t_5+t_5 - t_{10} = 1 - \frac1\beta \approx 0.38 \quad\text{but}\quad t_5+ t_5 - t_{11} = 1-\frac2\beta \approx - 0.24\,,
$$
which means that $\big|t_5+t_5 - (t_5\oplus t_5)\big| > \big|t_5+t_5 - t_{11}\big|$.
\end{ex}

The following example shows that also for $\beta$ root of $x^2-mx+1$, $m\geq 3$, the `sum' $t_j\oplus t_k$ is not always the closest $(-\beta)$-integer to $t_j+t_k$.

\begin{ex}
  Using the morphism $\overleftarrow{\varphi}_{-\beta}$ from Proposition~\ref{p:antimorfCP2}, one can list several letters of $u_{-\beta}$ around the delimiter
  marking the origin, namely
  $$
  u_{-\beta}=\cdots 1(0^{m-2}1)^{m-2}|\,0^{m-2}10^{m-3}1(0^{m-2}1)^{m-2}0^{m-2}10^{m-3}1(0^{m-2}1)^{m-2}0^{m-3}1(0^{m-2}1)^{m-2}\cdots\,.
  $$
Taking into account that $\Delta_0^-=1$ and $\Delta_1^-=2-\frac{n}{\beta}$, we derive the following.
 \begin{enumerate}
   \item  For the case $m=3$ we have $t_6=3\Delta_0^-+3\Delta_1^-=9-\frac{3}{\beta}$ and $t_6+t_6=18-\frac6\beta.$ The closest $(-\beta)$-integer to $t_6+t_6$ is $t_{11}=4\Delta_0^-+7\Delta_1^-=18-\frac7\beta$ instead of $t_{12}=t_{11}+\Delta_1^-=20-\frac7\beta.$
   \item  When $m\geq 4$ we have $t_2=2\Delta_0^-=2$ and $t_{m-2}=(m-2)\Delta_0^-=m-2$ and hence $t_2+t_{m-2}=m$
   while $t_m=(m-1)\Delta_0^-+\Delta_1^-=m+1-\frac1\beta.$ The closest $(-\beta)$-integer to $b_2+b_{m-2}$ is $b_{m-1}=m-\frac1\beta$.
   \end{enumerate}
\end{ex}

There is an interesting consequence of Proposition~\ref{p:soucet+}. The relation $-t_j=t_{-j}-\xi$ implies
$$
-\Big(t_j-\frac\xi2\Big)=t_{-j}-\frac\xi2\qquad\text{ for all }\ j\in\Z\setminus\{0\}\,.
$$
The geometrical interpretation of this fact is following.

\begin{coro}
  Let $\beta>1$ be the root of $x^2-mx-1,m\geq1.$ Then the set $\Z_{-\beta}\setminus\{0\}$ is symmetrical with respect to $\xi=\Delta_0^--\Delta_1^-.$
\end{coro}

The fact that the compatibility of the operations $\oplus$ and addition in $\R$ does not always hold is illustrated in the next example.
\begin{ex}
Let $\beta>1$ be the minimal Pisot number, zero of the polynomial $x^3-x+1$. Let us show that in this case the group operation $\oplus$
defined on $\beta$-integers by $t_j\oplus t_k=t_{j+k}$ is not compatible with addition in $\R$.

The R\'enyi expansion of 1 is equal to $d_\beta(1)=10001$, and thus the distances between consecutive $\beta$-integers take the values
$$
\Delta_0^+=1,\quad \Delta_1^+=\beta^{-4},\quad \Delta_2^+=\beta^{-3},\quad \Delta_3^+=\beta^{-2},\quad \Delta_4^+=\beta^{-1}.
$$
These values are not linearly independent over $\Q$. In particular, we have $\Delta_0^+=\Delta_2^++\Delta_3^+$.

The canonical substitution $\varphi_\beta$ for $\beta$ is of the form $0\mapsto 01$, $1\mapsto 2$, $2\mapsto 3$, $3\mapsto 4$, $4\mapsto 0$.
The infinite word $u_\beta$ coding the set $\Z_\beta$ is the fixed point of $\varphi_\beta$, namely
$$
\lim_{n\to\infty} \varphi^{n}(0) = 012340010120123012340123\cdots
$$
The first few $\beta$ integers can be written as
$$
\begin{aligned}
t_0&=0,\quad t_1=\Delta_0^+=1,\quad t_2 = \Delta_0^++\Delta_1^+ = 1+\frac1{\beta^4},\quad \\ t_3&=\Delta_0^++\Delta_1^++\Delta_2^+=1+\frac1{\beta^4}+\frac1{\beta^3},\quad\\
t_4&=\Delta_0^++\Delta_1^++\Delta_2^++\Delta_3^+=1+\frac1{\beta^4}+\frac1{\beta^3}+\frac1{\beta^2}\,\quad\dots
\end{aligned}
$$
We have $t_1\oplus t_2=t_3$, while $t_1+t_2 = 2\Delta_0^++\Delta_1^+ = \Delta_0^++\Delta_1^++\Delta_2^++\Delta_3^+ = t_4$.

\begin{figure}[ht]
{
\begin{center}
\setlength{\unitlength}{6pt}
\linethickness{0.7pt}
\begin{picture}(130,9)
\put(5,3){\line(1,0){62}}
\put(10,2.3){\line(0,1){1.4}}
\put(20,2.3){\line(0,1){1.4}}
\put(23.2,2.3){\line(0,1){1.4}}
\put(27.5,2.3){\line(0,1){1.4}}
\put(33.2,2.3){\line(0,1){1.4}}
\put(40.7,2.3){\line(0,1){1.4}}
\put(50.7,2.3){\line(0,1){1.4}}
\put(60.7,2.3){\line(0,1){1.4}}
\put(63.9,2.3){\line(0,1){1.4}}
%
%
\put(14.5,4){$\Delta_0^+$}
\put(20.8,4){$\Delta_1^+$}
\put(24.5,4){$\Delta_2^+$}
\put(29.5,4){$\Delta_3^+$}
\put(36,4){$\Delta_4^+$}
\put(45,4){$\Delta_0^+$}
\put(55,4){$\Delta_0^+$}
\put(61.3,4){$\Delta_1^+$}
\put(23.3,5.5){$\overbrace{\hspace*{59pt}}^{\text{\normalsize $t_1$}}$}
\put(9.7,5){$0$}
\put(9,0){$t_0$}
\put(19,0){$t_1$}
\put(23,0){$t_2$}
\put(27,0){$t_3$}
\put(33,0){$t_4$}
\put(40,0){$t_5$}
\put(50,0){$t_6$}
\put(60,0){$t_7$}
\put(63.5,0){$t_8$}
\put(66,0){$\cdots$}
\end{picture}
\end{center}
}
\caption{Operation $\oplus$ in $\Z_\beta$ is not compatible with addition in $\R$ for minimal Pisot number $\beta$.}
\label{f:minimalPisot}
\end{figure}
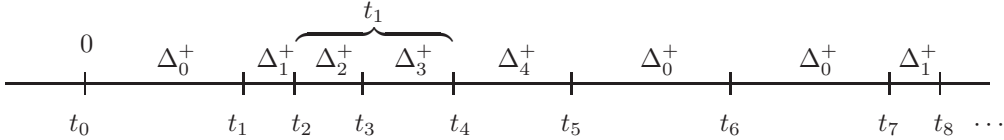
\end{ex}

The last example illustrates that compatibility does not enforce that the distance $|t_j+t_k-t_{j+k}|$ is bounded.
\begin{ex}
Let $\beta>1$ be the zero of the irreducible polynomial $x^6-x^5-1$. Such $\beta$ has conjugates outside of the unit circle,
so it is not a Pisot number. However, it is a Parry number. The R\'enyi expansion of 1 is $d_\beta(1)=100001$.
The distances between consecutive $\beta$-integers take the values
$$
\Delta_0^+=1,\quad \Delta_1^+=\beta^{-5},\quad \Delta_2^+=\beta^{-4},\quad \Delta_3^+=\beta^{-3},\quad \Delta_4^+=\beta^{-2},\quad \Delta_5^+=\beta^{-1},
$$
which are now linearly independent over $\Q$. By Theorem~\ref{t:compat}, the operation $\oplus$ defined on $\Z_\beta$ by $t_j\oplus t_k = t_{j+k}$ is compatible with addition in $\R$. However, the infinite word $u_\beta$ is fixed by a non-Pisot substitution $\varphi_\beta$.
The second largest eigenvalue $\lambda_2$ of its incidence matrix is in modulus greater than 1, and thus by~\cite{adamczewski}, the infinite word $u_\beta$
has unbounded balances. In particular, it means that $t_{j}+ t_{k}$ can be arbitrarily far from the ordinary sum $t_{j}\oplus t_{k}$. Nevertheless,
by Theorem~\ref{t:compat}, whenever $t_{j}+ t_{k}$ is a $\beta$-integer, it coincides with $t_{j}\oplus t_{k}$. Moreover, since the word $u_\beta$ is recurrent, i.e. each factor of $u_\beta$ (and the prefix of $u_\beta$ in particular) occurs infinitely many times, one can find arbitrarily large indices $j,k$ such that $t_j+t_k=t_{j+k}$.
\end{ex}

\section*{Acknowledgements}

\small
We would like to thank \v St\v ep\'an Starosta for very useful discussions. We acknowledge financial support by the Czech Science
Foundation grant 201/09/0584. The work was also partially supported by the CTU
student grant SGS11/085/OHK4/1T/14.

%
%


\end{document}